\documentclass[11pt]{amsart}
\usepackage{amsmath}
\usepackage{amssymb}
\usepackage{amscd}
\usepackage{amsthm}

\usepackage[all]{xy}

\usepackage{graphicx}
\usepackage{pstricks}
\usepackage{pst-plot}

\newtheorem{thm}{Theorem}[section]
\newtheorem{lem}[thm]{Lemma}
\newtheorem{cor}[thm]{Corollary}
\newtheorem{conj}[thm]{Conjecture}
\newtheorem{prop}[thm]{Proposition}

\newtheorem{claim}[thm]{Claim}

\theoremstyle{remark}

\newtheorem{rem}[thm]{Remark}

\theoremstyle{definition}
\newtheorem{defn}[thm]{Definition}

\numberwithin{equation}{section}

\allowdisplaybreaks[4]

\DeclareMathOperator{\Ext}{Ext}
\DeclareMathOperator{\Pic}{Pic}

\DeclareMathOperator{\ks}{ks}
\DeclareMathOperator{\Aut}{Aut}
\DeclareMathOperator{\cl}{cl}
\DeclareMathOperator{\CH}{CH}
\DeclareMathOperator{\im}{Im}

\begin{document}

\vfuzz0.5pc
\hfuzz0.5pc 

\newcommand{\claimref}[1]{Claim \ref{#1}}
\newcommand{\thmref}[1]{Theorem \ref{#1}}
\newcommand{\propref}[1]{Proposition \ref{#1}}
\newcommand{\lemref}[1]{Lemma \ref{#1}}
\newcommand{\coref}[1]{Corollary \ref{#1}}
\newcommand{\remref}[1]{Remark \ref{#1}}
\newcommand{\conjref}[1]{Conjecture \ref{#1}}
\newcommand{\questionref}[1]{Question \ref{#1}}
\newcommand{\defnref}[1]{Definition \ref{#1}}
\newcommand{\secref}[1]{Sec. \ref{#1}}
\newcommand{\ssecref}[1]{\ref{#1}}
\newcommand{\sssecref}[1]{\ref{#1}}
\newcommand{\figref}[1]{Figure \ref{#1}}

\def \red{{\mathrm{red}}}
\def \tors{{\mathrm{tors}}}
\def \EQ{\Leftrightarrow}

\def \mapright#1{\smash{\mathop{\longrightarrow}\limits^{#1}}}
\def \mapleft#1{\smash{\mathop{\longleftarrow}\limits^{#1}}}
\def \mapdown#1{\Big\downarrow\rlap{$\vcenter{\hbox{$\scriptstyle#1$}}$}}
\def \smapdown#1{\downarrow\rlap{$\vcenter{\hbox{$\scriptstyle#1$}}$}}
\def \A{{\mathbb A}}
\def \I{{\mathcal I}}
\def \J{{\mathcal J}}
\def \CO{{\mathcal O}}
\def \C{{\mathcal C}}
\def \BC{{\mathbb C}}
\def \BQ{{\mathbb Q}}
\def \m{{\mathcal M}}
\def \H{{\mathcal H}}
\def \BH{{\mathbb H}}
\def \S{{\mathcal S}}
\def \Z{{\mathcal Z}}
\def \BZ{{\mathbb Z}}
\def \W{{\mathcal W}}
\def \Y{{\mathcal Y}}
\def \T{{\mathcal T}}
\def \P{{\mathbb P}}
\def \CP{{\mathcal P}}
\def \G{{\mathbb G}}
\def \F{{\mathbb F}}
\def \BR{{\mathbb R}}
\def \D{{\mathcal D}}
\def \L{{\mathcal L}}
\def \f{{\mathcal F}}
\def \E{{\mathcal E}}
\def \BN{{\mathbb N}}
\def \N{{\mathcal N}}

\def \closure#1{\overline{#1}}
\def \EQ{\Leftrightarrow}
\def \imply{\Rightarrow}
\def \isom{\cong}
\def \embed{\hookrightarrow}
\def \tensor{\mathop{\otimes}}
\def \wt#1{{\widetilde{#1}}}

\def \ext{\mathop{{{\mathcal E}xt}}\nolimits}

\title[Indecomposable $K_{1}$ and Hodge-${\mathcal D}$ for $K3$ and Abelian]
{Indecomposable $K_{1}$ and the Hodge-${\mathcal 
D}$-conjecture for $K3$ and Abelian Surfaces}

\author{Xi Chen}

\address{632 Central Academic Building\\ 
University of Alberta\\
Edmonton, Alberta T6G 2G1, CANADA}
\email{xichen@math.ualberta.ca}

\author{James D. Lewis}

\address{632 Central Academic Building\\
University of Alberta\\
Edmonton, Alberta T6G 2G1, CANADA}
\email{lewisjd@gpu.srv.ualberta.ca}

\thanks{First author partially supported
by a startup grant from the University of Alberta. 
Second author partially supported
by a grant from the Natural Sciences and Engineering 
Research Council of Canada}
\subjclass{14C25, 14C30, 14C35}

\keywords{$K3$ surface, 
Abel-Jacobi map, regulator, Deligne cohomology, Chow group}

\date{Oct. 27, 2002}

\begin{abstract}
Let $X$ be a projective algebraic manifold, 
and $\text{CH}^{k}(X,1)$ the higher Chow group, with
corresponding real regulator 
$\text{r}_{k,1}\otimes {\BR}: \text{CH}^k(X, 1)\otimes {\BR}
\to H_{\mathcal D}^{2k-1}(X,{\BR}(k))$. If $X$ is a general $K3$
surface or Abelian surface, and $k=2$, we prove 
the Hodge-${\mathcal D}$-conjecture, i.e. the surjectivity
of $\text{r}_{2,1}\otimes {\BR}$. Since the 
Hodge-${\mathcal D}$-conjecture is not true for general
surfaces in $\P^{3}$ of degree $\geq 5$, the results in this
paper provide an effective bound for when this conjecture
is true. We then apply these results to the space of 
indecomposables ${\rm CH}_{\rm ind}^{k}(X,1;{\BQ})$, specifically
by proving that ${\rm Level}\big({\rm CH}_{\rm ind}^{k}(X,1;{\BQ})
\big) \geq k-2$ where $X$ is a general $k$-fold product 
of elliptic curves.  This leads to 
a hard generalization of Mumford's famous theorem
on the kernel of the Albanese map on the Chow group
of zero-cycles on a surface of positive genus.
\end{abstract}

\maketitle

\section{Statement of results}\label{SEC001}

Let $X$ be a  projective algebraic manifold. This
paper concerns the maps, called regulators, from $K_{1}$
of $X$ to real Deligne cohomology. More specifically,
in terms of Bloch's higher Chow groups $\text{CH}^{k}(X,m)$ \cite{Blo1},
we are interested in the case $m=1$ and the map
$$
\text{r}_{k,1}: \text{CH}^k(X, 1)\to H_{\mathcal D}^{2k-1}(X,{\BR}(k))
$$
where ${\BR}(k) = {\BR}(2\pi\sqrt{-1})^k$
and
\[
H_{\mathcal D}^{2k-1}(X,{\BR}(k)) \simeq H^{k-1,k-1}(X,{\BR}(k-1))
\]
is Deligne cohomology. Beilinson's 
Hodge-${\mathcal D}$-conjecture for real varieties would imply
that
$$
\text{r}_{k,1}\otimes {\BR}: \text{CH}^k(X, 1)\otimes {\BR}
\to H_{\mathcal D}^{2k-1}(X,{\BR}(k))
$$ 
is surjective (see \cite{Ja}). That conjecture is now known
to be false using the works of \cite{No} and \cite{G-S} (see \cite{MS1});
although the corresponding conjecture for smooth projective
varieties defined over number fields is still open. For
example if $S\subset {\P}^{3}$ is a transcendentally
generic surface of degree $d\geq 5$, then the image
$\text{cl}_{2,1}\big(\text{CH}^{2}(S,1)\big) \subset
H^{3}_{\mathcal D}(S,{\BR}(2))$ is contained
in the image of $H^{3}_{\mathcal D}({\P}^{3},{\BR}(2))
\to H^{3}_{\mathcal D}(S,{\BR}(2))$. The Hodge-${\mathcal D}$-conjecture
is trivially true for all smooth $S \subset {\P}^{3}$
of degree $d\leq 3$. This is because $H^{2}(S)$ is generated
by algebraic cycles if $d\leq 3$. When $d=4$, $S$ is a $K3$ surface,
and one can ask about the status of the Hodge-${\mathcal D}$-conjecture
here, and more generally for all $K3$ surfaces. In this paper
we settle this question for general $K3$ surfaces $S$, where general
means in the sense of the real analytic Zariski topology (see
\secref{SEC002}). In other words $H^{1,1}(S,{\BR}(1)) \simeq \
H^{1,1}(S,{\BR}(1))^{\vee}$ is generated by currents of the form
$$
\omega\in H^{1,1}(S,{\BR}(1)) \mapsto \frac{1}{2\pi\sqrt{-1}}\sum_{j}
\int_{D_{j}}\omega\log|f_{j}|,
$$
where $f_{j}\in \BC(D_{j})^{\times}$ and
$\sum_{j}\text{div}(f_{j}) = 0$.

\bigskip

More specifically, we prove the following:

\begin{thm}\label{THM100}
{\rm (i)} {\it The Hodge-${\mathcal D}$-conjecture holds
for general $K3$ surfaces in the real analytic Zariski topology.

\medskip

{\rm (ii)} The Hodge-${\mathcal D}$-conjecture holds
for general Abelian surfaces in the real analytic Zariski topology,
hence for general Kummer surfaces.

\medskip

{\rm (iii)} The Hodge-${\mathcal D}$-conjecture holds
for general products $E_{1}\times E_{2}$ of elliptic
curves in the real analytic Zariski topology, hence
for general ``special Kummer surfaces'' (see \secref{SEC002}).}
\end{thm}

There is ample evidence why such a statement should be true. 
There are the works of A. Collino, S. M\"uller-Stach, C. Voisin, 
{\it et al} on nontrivial regulator calculations on $K3$
surfaces. See for example \cite{Co2}, \cite{MS1}
and the references cited there. Further, the second author proved a 
twisted version of the conjecture with twisted
higher Chow groups and
regulators \cite{Lw1}. 
There are two key ingredients which make a proof
for a general $K3$ surface $S$
possible. First of all, there are plenty of nontrivial higher Chow
cycles on $S$ constructed out of rational curves. Take two rational
curves $D_1$ and $D_2$ and two distinct points $p_1, p_2\in D_1\cap
D_2$. Choose $f_i\in \BC(D_i)^\times$ such that
\begin{equation}\label{E005}
(f_1) = p_1 - p_2  \text{ and } (f_2) = p_2 - p_1
\end{equation}
and then $(f_1, D_1) + (f_2, D_2)$ defines
a class in $\text{CH}^2(S, 1)$. More generally, take
$n$ rational curves $D_1, D_2, \ldots, D_n$ and pick $n$ distinct points
$p_1, p_2,\ldots, p_n$ such that $p_i\in D_i\cap D_{i+1}$
(let $D_{n+j} = D_j$ and $p_{n+j} = p_j$).
Choose $f_i\in \BC(D_i)^\times$ such that
$$
(f_i) = p_i - p_{i+1}
$$
and then $\sum_{i=1}^n (f_i, D_i)\in \text{CH}^2(S, 1)$.

It is well known that there are rational curves on $S$.
This statement was made more precise in \cite{C1}, that there
are rational curves in every linear series on $S$. Indeed, we think
the following is true:

\begin{conj}
The rational curves on $S$ are dense. That is, let
$\Sigma = \cup D_\alpha$ be
the union of all rational curves on $S$. Then the closure
$\overline{\Sigma}$ of $\Sigma$ under the analytic topology is $S$.
\end{conj}

The motivation for this conjecture is based on the analogy
between rational curves on $K3$ surfaces with
the density of torsion points on elliptic curves. These
torsion points were instrumental in constructing nontrivial
$K_{2}$ classes on general elliptic curves (see \cite{Blo2}),
and correspondingly, the rational curves on a general $K3$ surface
play a role here in constructing nontrivial $K_{1}$ classes.

\bigskip

Actually, there are enough rational curves even only in the primitive
class by the counting of \cite{B-L} and \cite{Y-Z} (a primitive class
is a divisor which is not the multiple of another divisor; there is
only one primitive class on a general algebraic K3, which has
self-intersection $2g - 2$ with $g$ called the genus of the K3).
For example, there
are $3200$ rational curves in the primitive class of a quartic $K3$,
which seem enough to produce $20$ generators of
$H^{1,1}(S, {\BR})$. It turns out in the end that we only need a
fraction of these $3200$ rational curves to realize the surjection of
$\text{cl}_{2,1}\otimes {\BR}$.

The second key fact is that Hodge-${\mathcal D}$-conjecture obviously holds on a
$K3$ surface with maximum Picard number $20$. This points us to a
way to prove the conjecture for general $K3$ surfaces by
degeneration. Actually, this is also the way in which the twisted
version of the conjecture was proven \cite{Lw1}.

\bigskip

In \secref{SEC005} we turn our attention to indecomposability of $K_{1}$
for special classes of varieties. There we introduce the
notion of Level, which measures the complexity of the Chow groups
$\text{CH}^{k}(X,m)$, as well as the quotient group
of indecomposables $\text{CH}^{k}_\text{ind}(X,m)$. By exploiting
the results in \thmref{THM100}, we arrive at the following:

\begin{thm}\label{THM101}
Let $X = E_{1}\times \cdots \times E_{k}$ 
be a general product of $k$ elliptic curves, i.e. in the countable
real analytic Zariski topology (see \secref{SEC002}), 
and ${\rm CH}_{\rm ind}^{k}(X,1;{\BQ})$ the
space of indecomposables. Then:

\medskip

{\rm (i)} ${\rm Level}\big({\rm CH}_{\rm ind}^{k}(X,1;{\BQ})\big) \geq k-2$.

\medskip

{\rm (ii)} In particular for $k\geq 3$, there are an 
uncountable number of indecomposables
in the kernel of the regulator 
$$
{\rm cl}_{k,1}: {\rm CH}_{\rm ind}^{k}(X,1;{\BQ})
\to H_{\mathcal D}^{2k-1}(X,{\BQ}(k)).
$$
\end{thm}

This paper has its origins in a discussion that the second
author had with the first, where it was suggested
by the second author
that the Hodge-${\mathcal D}$-conjecture for $K3$ surfaces should
be true, based on the abundance of rational curves, and that
a proof should involve degeneration to a $K3$ surface of
maximum Picard number. The second author is indeed very grateful to the
first author for supplying the complete degeneration argument
in \secref{SEC008}, without which this paper would not have evolved in
its present form.

One consequence of the results of this paper
is a significant generalization of Theorems 1 and 2 
on page 544 of \cite{GL1} (and corresponding statements in
\cite{GL2}). Not only do we present correct proofs of
these theorems, the results in this paper are much
deeper. The second author is grateful to Morihiko Saito
for pointing out the errors in the degeneration argument
in \cite{GL1}, where  the cycles constructed in Theorem
2.4 and Proposition 3.3 of \cite{GL1} are regulator decomposable,
contrary to the claims there. The problem in \cite{GL1}
has to do with the presence of singularities of a real $2$-form 
after degeneration to a singular fiber. Thus this paper
can also be seen as providing the  correct proofs to the main 
results in \cite{GL1} (and \cite{GL2}).

\section{Notation}\label{SEC002}

Throughout this paper, $X$ is assumed to
be a projective algebraic manifold of dimension $n$. If
${\bf A} \subset {\BR}$ is a subring, we put 
${\bf A}(k) = {\bf A}(2\pi\sqrt{-1})^k$. For the higher
Chow groups $\text{CH}^{k}(X,m)$ introduced in \secref{SEC003}, and for
${\bf A}$ above, we denote $\text{CH}^{k}(X,m)\otimes {\bf A}$
by $\text{CH}^{k}(X,m;{\bf A})$.

The use of the terminology ``general $X$'' in this paper
will have two possible meanings. Firstly, for a variety $Y$,
a real analytic Zariski open set $U$ in $Y$ will be the
complement of a real analytic subvariety of $Y$. If $Y$
parameterizes a family $\{X_{t}\}_{t\in Y}$ of projective 
algebraic manifolds, then a general such $X = X_{t}$ 
in the real analytic Zariski topology means that $t\in U$,
for some real analytic Zariski open set $U$ in $Y$. 
Secondly, a countable
real analytic open set $U_{c}$ of $Y$ is the complement
of a countable union of real analytic subvarieties of $Y$.
We say that $X = X_{t}$ is general in the countable real
analytic Zariski topology, if $t\in U_{c}$.

Let $A$ be an Abelian surface, and $Y = A/\pm 1$ its corresponding
Kummer counterpart \cite[p. 550]{S-S}. Following \cite{S-S}, we say
that $Y$ is special if $A$ reducible, i.e. contains an elliptic
curve.

\section{Deligne cohomology and higher Chow groups}\label{SEC003}

{\it (a) Deligne cohomology.\/}
Let 
$$
\Omega_{X}^{\bullet} := {\mathcal O}_{X}\ {\buildrel d\over\to}\ 
\Omega_{X}^{1}\ {\buildrel d\over\to}\ \Omega_{X}^{2}\
{\buildrel d\over\to}\ 
\cdots \ {\buildrel d\over\to}\  \Omega_{X}^{n}\to 0,
$$
be the holomorphic de Rham complex.
The complex $\Omega_{X}^{\bullet}$ is filtered by the Hodge filtration
$$
F^{p}\Omega_{X}^{\bullet} := 0\to \Omega_{X}^{p}\ {\buildrel d\over\to}\ 
\Omega_{X}^{1}\ {\buildrel d\over\to}\ \Omega_{X}^{2}\ {\buildrel d\over\to}\ 
\cdots \ {\buildrel d\over\to}\  \Omega_{X}^{n}\to 0.
$$
The truncated complex $\Omega_{X}^{<p}$ is defined by the short
exact sequence
$$
0\to F^{p}\Omega_{X}^{\bullet}\to \Omega_{X}^{\bullet}\to
\Omega_{X}^{<p}\to 0.
$$
We recall that the Hodge spectral sequence
$$
H^{q}(X,\Omega_{X}^{p})\ \Rightarrow\ H_\text{DR}^{p+q}(X,\BC)
$$
degenerates at $E_{1}$. Thus
$$
F^{p}H_\text{DR}^{k}(X,\BC) := 
{\BH}^k(F^{p}\Omega_{X}^{\bullet})\hookrightarrow
{\BH}^{k}(\Omega_{X}^{\bullet})
$$
is injective for all $k$. 
Therefore if we put
$\Omega_{X}^{<p} := \Omega_{X}^{\bullet}/F^{p}\Omega_{X}^{\bullet}$,
then 
$$
{\BH}^{k}(\Omega_{X}^{<p}) = 
\frac{H^{k}_\text{DR}(X,\BC)}{F^{p}H^{k}_\text{DR}(X,\BC)}.
$$

For a subring ${\bf A}\subset {\BR}$, we introduce the Deligne complex
$$
{\bf A}_{\mathcal D}(k):\quad  {\bf A}(k) \to \underbrace{{\mathcal O}_X \to 
\Omega^1_X\to\cdots \to \Omega^{k-1}_X}_{\text{call this}\ 
\Omega_X^{\bullet < k}}.
$$

\begin{defn}\label{DEF300}
Deligne cohomology is given by the hypercohomology:
$$
H_{\mathcal D}^{i}(X,{\bf A}(k)) ={\BH}^{i}({\bf A}_{\mathcal D}(k)).
$$
\end{defn}

Applying ${\BH}^\bullet(-)$ to the short exact sequence:
$$
0\to \Omega_X^{\bullet < k}[-1]\to {\bf A}_{\mathcal D}(k) \to {\bf A}(k)\to 0,
$$
yields the short exact sequence:
\[
\begin{split}
0 & \to \frac{H^{i-1}(X,\BC)}{H^{i-1}(X,{\bf A}(k)) +
F^kH^{i-1}(X,\BC)}\\
& \to H_{\mathcal D}^{i}(X,{\bf A}(k)) 
\to H^i(X,{\bf A}(k))\cap F^kH^i(X,\BC) \to 0.
\end{split}
\]
The cases of interest are ${\bf A} = {\BZ}$, ${\BQ}$ 
and ${\BR}$, for $i=2k-1$.
For example, if ${\bf A} = {\BZ}$ and $i=2k-1$, then we arrive at the
short exact sequence
\[
\begin{split}
0 &\to \frac{H^{2k-2}(X,\BC)}{F^{k}H^{2k-2}(X,\BC)
+ H^{2k-2}(X,{\BZ}(k))}\\
& \to H_{\mathcal D}^{2k-1}(X,{\BZ}(k))
\to H^{2k-1}(X,{\BZ}(k))_\text{tor} \to 0.
\end{split}
\]

Next, if ${\bf A} = {\BQ}$ and $i= 2k-1$, then from Hodge
theory, 
$$
H^{i}(X,{\BQ}(k))\cap F^kH^{i}(X,\BC) = 0.
$$
Thus we arrive at the isomorphism
$$
H^{2k-1}_{\mathcal D}(X,{\BQ}(k)) \simeq
\frac{H^{2k-2}(X,\BC)}{F^{k}H^{2k-2}(X,\BC) + H^{2k-2}(X,{\BQ}(k))}.
$$
If ${\bf A} = {\BR}$ and if we set 
$$
\pi_{k-1} : \BC = {\BR}(k)\oplus {\BR}(k-1) \to {\BR}(k-1)
$$
to be the projection, then we have the isomorphisms:
\[
\begin{split}
H^{2k-m}_{\mathcal D}(X,{\BR}(k)) &\simeq
\frac{H^{2k-m-1}(X,\BC)}{F^kH^{2k-m-1}(X,\BC) +
H^{2k-m-1}(X,{\BR}(k))}\\
&
\xrightarrow[\simeq]{\pi_{k-1}}
\frac{H^{2k-m-1}\big(X,{\BR}(k-1)\big)}{\pi_{k-1}
\big(F^{k}H^{2k-m-1}(X,\BC)\big)}.
\end{split}
\]
For example if ${\bf A} = {\BR}$ and $i=2k-1$, we have
\[
\begin{split}
H^{2k-1}_{\mathcal D}(X,{\BR}(k)) & \simeq
\frac{H^{2k-2}(X,\BC)}{F^k H^{2k-2}(X,\BC) + H^{2k-2}(X,{\BR}(k))}\\
& \xrightarrow[\simeq]{\pi_{k-1}} H^{k-1,k-1}(X,{\BR})\otimes
{\BR}(k-1)\\
&=: H^{k-1,k-1}(X,{\BR}(k-1))\\
&\simeq \biggl\{H^{n-k+1,n-k+1}(X, {\BR}(n-k+1))\biggr\}^\vee.
\end{split}
\]

{\it (b) Higher Chow groups.\/}
The higher Chow groups $\text{CH}^{k}(X,m)$ were invented
by Bloch \cite{Blo1} (and independently by S.
Landsberg). They are defined for all $k,\ m \geq 0$; moreover
in \cite{Blo1} there is proven a Riemann-Roch theorem, namely that the
Chern character map
$$
\text{ch} : \text{CH}^{\bullet}(X,m;{\BQ}) \
{\buildrel\sim\over \to}\ K_{m}(X)\otimes {\BQ},
$$
is an isomorphism, generalizing the Grothendieck Riemann-Roch
theorem in the case $m=0$. Since in this paper, we are only interested
in the case $m=1$, we provide an abridged alternative 
version of the definition of Bloch's higher Chow groups, for
the cases $0\leq m\leq 2$, that we aquire using a Gersten resolution. Let
${\mathcal K}_{k,X}$ be the sheaf of $K$-groups on $X$, i.e. where
where ${\mathcal K}_{k,X}$ is the sheaf associated to the
presheaf in the Zariski topology, 
$$
U\subset X \mapsto K_{k}\big(\Gamma(U,{\mathcal O}_{X})\big).
$$

\begin{defn}\label{DEF301}
(See \cite{MS1}) For $0\leq m\leq 2$, define
$$
\text{CH}^{k}(X,m) = H_\text{Zar}^{k-m}(X,{\mathcal K}_{k,X}).
$$
\end{defn}

The best way to interpret the RHS in \defnref{DEF301} is
via the Gersten resolution proven by Bloch for $k=2$
and by Quillen for general $k$. This is the flasque 
resolution of ${\mathcal K}_{k,X}$ given by
\[
\begin{split}
0 &\to {\mathcal K}_{k,X}\to K_{k}(\BC(X)) \to \bigoplus_{\text{
cd}_{X}Z=1}i_{Z,\ast}K_{k-1}(\BC(Z))\to\cdots\\
& \to \bigoplus_{\text{cd}_{X}Z=k-2}i_{Z,\ast}K_{2}(\BC(Z))\to
\bigoplus_{\text{cd}_{X}Z=k-1}i_{Z,\ast}K_{1}(\BC(Z))\\
&\to\bigoplus_{\text{cd}_{X}Z=k}i_{Z,\ast}K_{0}(\BC(Z))\to 0
\end{split}
\]

We recall that $K_{0}(\BC(Z))\simeq {\BZ}$,
$K_{1}(\BC(Z)) = \BC(Z)^{\times}$, and that
$K_{2}(\BC(Z))$ is generated by symbols. Taking global 
sections leads to a complex, whose last three terms are
$$
\bigoplus_{\text{cd}_{X}Z=k-2}K_{2}(\BC(Z))\
{\buildrel T\over\to}\ \bigoplus_{\text{cd}_{X}Z=k-1}
\BC(Z)^{\times}\ {\buildrel \text{div}\over\to}\ z^{k}(X),
$$
where $T$ is the Tame symbol, and div is the divisor map, and
$z^{k}(X)$ is the free abelian group generated by subvarieties
of codimension $k$ in $X$. Thus for example,
\[
\text{CH}^{k}(X,0) = \frac{z^{k}(X)}{\left\{
\begin{matrix}
\text{subgroup of}\\
\text{principal divisors}
\end{matrix}
\right\}} 
= \text{CH}^{k}(X),
\]
and
\[
\text{CH}^{k}(X,1) = \frac{\left\{
\sum_{j}(f_{j},Z_{j})\ :\
\begin{matrix}
\text{cd}_{X}Z_{j} = k-1,\ f_{j}\in \BC(Z_{j})^{\times}\\
\sum_{j} \text{div}(f_{j}) = 0
\end{matrix}
\right\}}{\text{Image(Tame symbol)}}.
\]

\section{A regulator}\label{SEC004}

In this section, we recall the definition of the regulator
$$
\text{c}_{k,1}: \text{CH}^{k}(X,1) \to H_{\mathcal D}^{2k-1}(X,
{\BZ}(k)),
$$
where we recall that  $H_{\mathcal D}^{2k-1}(X,{\BZ}(k))$
fits in the short exact sequence:
\[
\begin{split}
0 & \to \frac{H^{2k-2}(X,\BC)}{F^{k}H^{2k-2}(X,\BC)
+ H^{2k-2}(X,{\BZ}(k))} \to H_{\mathcal D}^{2k-1}(X,{\BZ}(k))\\
& \to H^{2k-1}(X,{\BZ}(k))_\text{tor} \to 0.
\end{split}
\]
We define 
$$
\text{CH}_{\hom}^{k}(X,1) = \ker \biggl(\text{CH}^{k}(X,1) \to
H^{2k-1}(X,{\BZ}(k))\biggr).
$$
Using the compatibility of Poincar\'e and Serre duality,
there is an induced map
$$
\underline{\text{cl}}_{k,1}: \text{CH}_{\hom}^{k}(X,1) 
\to 
\frac{\{F^{n-k+1}H^{2n-2k+2}(X,{\BC})\}^{\vee}}{H_{2n-2k+2}(X,{\BZ}(n-k))}.
$$
The formula we use for
$\underline{\text{cl}}_{k,1}$ is this (see \cite{Lev}): Consider
$$
\xi := \sum_{j}(f_{j}, D_{j}), \quad \sum_{j}\text{div}(f_{j}) = 0,
$$
where $\text{codim}_{X}D_{j} = k-1$ and $f_{j}\in \BC(D_{j})^{\times}$.
Choose a branch of the $\log$ function on $\BC\backslash
[0,\infty)$, and put $\gamma_{j} = f_{j}^{-1}[0,\infty]$, and
$\gamma = \sum_{j}\gamma_{j}$. Then $\partial \gamma = 0$ on
$X$, and our assumption that $\xi \in \text{CH}^{k}_{\hom}(X,1)$
means that $\gamma = \partial \zeta$ is a boundary. 
For $\omega\in F^{n-k+1}H^{2n-2k+2}(X,\BC)$,
$$
\underline{\text{cl}}_{k,1}(\xi)(\omega) =
\frac{1}{(2\pi\sqrt{-1})^{n-k+1}}
\biggl(\sum_{j}\int_{D_{j}\backslash 
\gamma_{j}}\omega\log(f_{j})\ +\
2\pi\sqrt{-1}\int_{\zeta}\omega\biggr).
$$
Now consider a smooth family of projective algebraic manifolds 
$Y := \coprod_{t\in \Delta}X_{t}$ over a disk $\Delta \subset \BC$.
We want to show that $\underline{\text{cl}}_{k,1}$ varies holomorphically
with respect to a family of $K_{1}$-cycles in $Y$ over $\Delta$. Assume
given an algebraic family of cycles
$\xi_{t} =  \sum_{j}(f_{j,t},D_{j,t})\in \text{CH}_{\hom}^{k}(X_{t},1)$, $t\in \Delta$. 
As in \cite{Gr}, one can choose $C^{\infty}$ differential forms
$\omega_{1},\ldots,\omega_{r}$ on $Y$ such that each $\omega_{i}$
is of Hodge type $(n,n-2k+2)+\cdots +(n-k+1,n-k+1)$, $d\omega_{i}\wedge
dt = 0$, and $\{\omega_{i}\big|_{X_{t}}\ :\ i=1,\ldots, r\}$ gives
a basis of $F^{n-k+1}H^{2n-2k+2}(X_{t},\BC)$ for all $t\in \Delta$.
Let $\omega$ be any linear combination of $\omega_{1},\ldots,\omega_{r}$.
This next result is probably well-known to experts, although
a proof does not seem to be written down in the literature. 
For the convenience to the reader, we present a proof here.

\begin{prop}\label{PROP400}
$\underline{\text{cl}}_{k,1}(\xi_{t})(\omega)$ varies
holomorphically in $t\in \Delta$.
\end{prop}

\begin{proof}
We base our proof on similar ideas in Appendix A
of \cite{Gr}. First, it is reasonably clear that
$\underline{\text{cl}}_{k,1}(\xi_{t})(\omega)$ varies continuously
in $t\in \Delta$. We use the criterion of holomorphicity via
Morera's theorem. Let $\Gamma \subset \Delta$ be a simple-closed
curve oriented counterclockwise. We need to show that
$\int_{\Gamma}\underline{\text{cl}}_{k,1}(\xi_{t})(\omega) dt = 0$.
This amounts to the calculation of
$$
\sum_{j}\int_{\Gamma}\biggl(\int_{D_{j,t}\backslash 
\gamma_{j,t}}\omega\log(f_{j,t})\biggr)dt\ +\
2\pi\sqrt{-1}\int_{\Gamma}\biggl(\int_{\zeta_{t}}\omega\biggr) dt.
$$
Note that $\Gamma = \partial M$ for some region $M\subset \Delta$.

\bigskip

(I) We first evaluate $2\pi\sqrt{-1}\int_{\Gamma}\big(
\int_{\zeta_{t}}\omega\big)  dt$: 
Put $\zeta_{\Gamma} = \coprod_{t\in \Gamma}\zeta_{t}$,
$\zeta_{M} = \coprod_{t\in M}\zeta_{t}$ and $\gamma_{M}
= \coprod_{t\in M}\gamma_{t}$. Then by Stokes' theorem
\[
2\pi\sqrt{-1}\int_{\Gamma}\biggl(\int_{\zeta_{t}}\omega\biggr)dt=
\int_{\zeta_{\Gamma}}\omega\wedge dt = 
-2\pi\sqrt{-1}\int_{\gamma_{M}}\omega\wedge dt.
\]

(II) Next we evaluate $\sum_{j}\int_{\Gamma}\big(\int_{D_{j,t}\backslash 
\gamma_{j,t}}\omega\log(f_{j,t})\big)dt$: Put
\[
D_{j,\Gamma} = 
\coprod_{t\in \Gamma}D_{j,t} \text{ and }
D_{j,M} = \coprod_{t\in M}D_{j,t}.
\]
Note that our assumptions on $\omega$, (and holomorphicity of
$f_{j,t}$, away from the pole sets) imply that 
$$
d\big(\log(f_{j,t})\omega\wedge dt\big) = \frac{df_{j,t}}{f_{j,t}}\wedge
\omega\wedge dt,
$$
and that by Hodge type considerations alone, 
$$
\biggl(\frac{df_{j,t}}{f_{j,t}}\wedge
\omega\wedge dt\biggr)\biggl|_{D_{j,M}} \ =\ 0.
$$
Thus by Stokes' theorem, we have
\[
\begin{split}
\sum_{j}\int_{\Gamma}\biggl(\int_{D_{j,t}\backslash 
\gamma_{j,t}}\omega\log(f_{j,t})\biggr)dt
&= \sum_{j}\int_{D_{j,\Gamma}} \log(f_{j,t})\omega\wedge dt\\
&= 2\pi\sqrt{-1}\int_{\gamma_{M}}\omega\wedge dt,
\end{split}
\]
where we use the fact that we pick up a $2\pi\sqrt{-1}$ ``period'' from
$\log(f_{j,t})$ as we cross $\gamma_{j,t}$.

\bigskip

Finally, the resulting terms from (I) and (II) cancel, which 
establishes the proposition.
\end{proof}

We note that the forms $\omega_{1},\ldots,\omega_{r}$ define
a holomorphic frame of the bundle 
$$
\coprod_{t\in \Delta}F^{n-k+1}H^{2n-2k+2}(X_{t},\BC)
$$
over $\Delta$. One can
define a frame of the $C^{\infty}$ ${\BR}(n-k+1)$-bundle 
$$
\coprod_{t\in \Delta}H^{n-k+1,n-k+1}(X_{t},{\BR}(n-k+1))
$$ 
in terms of a linear combinations of the forms 
$\omega_{1},\ldots,\omega_{r}$.
The coefficient functions of these linear combinations will
be complex-valued combinations of real analytic functions. 
If we put $X = X_{t}$, then this
follows from the fact that forms in $H^{n-k+1,n-k+1}(X,\BC)$
are by duality, precisely subspace of the forms in 
$F^{n-k+1}H^{2n-2k+2}(X,\BC)$
killed by the wedge products of forms in 
$\overline{F^{k}H^{2k-2}(X,\BC)}$ (see \remref{REM403}).
Further, for any form $\omega\in H^{n-k+1,n-k+1}(X,{\BR}(n-k+1))$,
\[
\begin{split}
&\quad\text{Re}\left(\frac{1}{(2\pi\sqrt{-1})^{n-k+1}}
\biggl(\sum_{j}\int_{D_{j}\backslash 
\gamma_{j}}\omega\log(f_{j})\ +\
2\pi\sqrt{-1}\int_{\zeta}\omega\biggr)\right)\\
&= \frac{1}{(2\pi\sqrt{-1})^{n-k+1}}\sum_{j}\int_{D_{j}}\omega\log|f_{j}|.
\end{split}
\]

We denote by $\text{r}_{k,1}$, the corresponding Beilinson real regulator
\begin{equation}\label{E401}
\begin{split}
\text{r}_{k,1} : \text{CH}^{k}(X,1) & \to H^{2k-1}(X,{\BR}(k))\\
&\simeq 
H^{n-k+1,n-k+1}(X,{\BR}(n-k+1))^{\vee}
\end{split}
\end{equation}
given by
$$
\xi = \sum_{j}(f_{j},D_{j})\mapsto \text{r}_{k,1}(\xi)(\omega)
= \frac{1}{(2\pi\sqrt{-1})^{n-k+1}}\sum_{j}\int_{D_{j}}\omega\log|f_{j}|.
$$

\begin{cor}\label{COR402}
The corresponding real regulator defines a real
analytic function in $t\in \Delta$. In particular, if we are given
a smooth family of projective varieties
${\mathcal X} := \coprod_{t\in S}X_{t}  \to S$,
over a smooth quasiprojective base $S$, and family of 
of $K_{1}$ classes $\xi_{t}\in \text{CH}^{k}(X_{t},1)$ with
nonvanishing real regulator value $\text{r}_{k,1}(\xi_{t_{0}})$, then
$\text{r}_{k,1}(\xi_{t}) \ne 0$ for $t$ in some nonempty real analytic Zariski
open subset of $S$.
\end{cor}

\begin{rem}\label{REM403}
Here
are some details on the existence of real analytic frames:
One can assume that the first
$\{\omega_{1},\ldots,\omega_{\ell}\}
\subset \{\omega_{1},\ldots,\omega_r\}$ restrict to a basis
for
\[
F^{n-k+2}H^{2n-2k+2}(X_{t},\BC).
\]
A corresponding
holomorphic frame gives $C^{\infty}$ forms 
$\eta_{1},\ldots,\eta_{\ell}$
restricting to a basis of $F^{k}H^{2k-2}(X_{t},\BC)$. One
can construct a $C^{\infty}$ frame restricting to a basis
of $H^{n-k+1,n-k+1}(X_{t},{\BR}(n-k+1))$ by finding
the general solution of the linear
system
\[
\langle a_{1}(t)\omega_{1}+\cdots + a_{r}(t)\omega_{r},
\overline{\eta}_{j}\rangle = 0
\text{ for }
j = 1,\ldots,\ell,
\]
together with the nonsingularity of the matrix $\big(\langle \omega_{i},
\overline{\eta}_{j}\rangle_{1\leq i,j\leq \ell}\big)$. After which, via the
projection 
$\pi_{n-k+1} : \BC\to {\BR}(n-k+1),$ we arrive
at the real analytic frame, twisted by ${\BR}(n-k+1)$ as
required.
\end{rem}

\section{Indecomposability}\label{SEC005}

In \cite{Lw2} we introduced the notion of the Level of
a Chow group (see \label{DEFN503} below). Similar notions appear
elsewhere, such as in \cite{Sa}. One should think of the Level
as an integral invariant measuring the complexity
of a given Chow group. Using this notion of Level, we will
see that $\text{CH}^k_\text{ind}(X,1;{\BQ})$ can be very
large (and uncountable); moreover even in the kernel of the regulator
map. Thus one arrives at an analogous result to Mumford's
famous theorem \cite{Md}. 

\bigskip

The setting is this. Recall Bloch's higher Chow group
$\text{CH}^k(X,m)$ \cite{Blo1}. As in the case 
$m \leq 1$, Bloch \cite{Blo3} (as well as Beilinson)
constructs a cycle class map
$$
\text{cl}_{k,m} : \text{CH}^k(X,m) \rightarrow H_{\mathcal D}^{2k-m}
(X,{\BZ}(k)).
$$
For $m \geq 1$, we have for example
$$
H_{\mathcal D}^{2k-m}(X,{\BQ}(k))
\simeq \frac{H^{2k-m-1}(X,\BC)}{F^kH^{2k-m-1}(X,\BC) +
H^{2k-m-1}(X,{\BQ}(k))}.
$$
One has products $\text{CH}^k(X,m) \cap \text{CH}^r(X,{\ell}) \subset
\text{CH}^{k+r}(X,m+\ell)$ compatible with the product structure
on Deligne cohomology.

\bigskip

We now assume $m \geq 1$.

\begin{defn}\label{DEF500}
(i) The subgroup of decomposables is given by
\[
\text{CH}_D^k(X,m) := \text{Image}\biggl((\BC^{\times})^{\otimes_{\BZ}m}
\otimes_{\BZ}\text{CH}^{k-m}(X,0) \ {\buildrel \cap \over \longrightarrow}\
\text{CH}^k(X,m)\biggr),
\]
where $\BC^{\times}$ is identified with $\text{CH}^1(X,1)$
via the identification 
$$
\text{CH}^1(X,1) = H^{0}_\text{Zar}(X,
{\mathcal K}_{1,X}) =  H^{0}_\text{Zar}(X,{\mathcal O}^{\times}_{X}).
$$

\medskip

(ii) The space of indecomposables is given by
$$
\text{CH}^k_\text{ind}(X,m;{\BQ}) := \text{CH}^k(X,m;{\BQ}) 
/ \text{CH}^k_D(X,m;{\BQ}).
$$
\end{defn}

\bigskip

\begin{rem}\label{REM501}
(i) There is an isomorphism
\[
\begin{matrix}
\text{cl}_{1,1} &: &\text{CH}^1(X,1)&{\buildrel \sim \over\rightarrow}&
H^1_{\mathcal D}(X,{\BZ}(1))&\simeq&H^0(X,\BC/{\BZ}(1))\\
&\\
&&\vert\vert&&&&\vert\vert\\
&\\
&&\ \BC^{\times}&&{\buildrel \text{Id}\over\rightarrow}
&&\ \BC^{\times}\\
\end{matrix}
\]

\bigskip

(ii) The product structure in Deligne cohomology implies
that
$$
H^1_{\mathcal D}(X,{\BZ}(1))\ \cup\ H^1_{\mathcal D}(X,{\BZ}(1))
= 0 \in H^2_{\mathcal D}(X,{\BZ}(2)).
$$
Therefore $\text{cl}_{k,m}(\text{CH}^k_D(X,m)) = 0$ for $m \geq 2$.

\bigskip

(iii) In the case $m=1$, we have
$$
\text{cl}_{k,1}(\text{CH}^k_D(X,1;{\BQ})) \subset \BC^{\times}_{\BQ}
\otimes H^{k-1,k-1}(X,{\BQ}(k-1)), \ \text{where} \ \BC^{\times}_{\BQ}
= \BC / {\BQ}(1),
$$
with equality $\Leftrightarrow$ the Hodge conjecture holds for $X$.

\bigskip

(iv) A rigidity result of Beilinson \cite{Bei1} implies that
the image  
\[
\text{cl}_{k,m}(\text{CH}^k(X,m;{\BQ}))
\]
is countable for
$m \geq 2$. A variant of this rigidity argument \cite{MS1}
shows that the induced map
$$
\underline{\text{cl}}_{k,1} :
\frac{\text{CH}^k(X,1;{\BQ})}{\text{CH}^k_D(X,1;{\BQ})}
\rightarrow \frac{H^{2k-1}_{\mathcal D}(X,{\BQ}(k))}{
\BC^{\times}_{\BQ} \otimes H^{k-1,k-1}(X,{\BQ}(k-1))},
$$
has countable image.
\end{rem}

Before stating our main results, we introduce some terminology.

$N^{\bullet}H^{\ast}(X,{\BQ}) =$ filtration by coniveau, with graded
piece
$$
\text{Gr}_N^{\ell}H^i(X,{\BQ}) = \frac{N^{\ell}H^i(X,{\BQ})}{
N^{\ell + 1}H^i(X,{\BQ})}.
$$

More explicitly,
\[
N^jH^i(X,{\BQ}) = \ker : H^i(X,{\BQ}) \to
\lim_{\buildrel \longrightarrow\over {
\begin{matrix}
{}^{Y\hookrightarrow X \
\text{closed}}\\
{}^{\text{codim}_XY \geq j}
\end{matrix}
}}
H^i(X-Y,{\BQ}).
\]

\begin{defn}\label{DEF503}
Let $G$ be a subgroup of $\text{CH}^k(X,m;{\BQ})$.
Then

\medskip

(i) $\text{Level}(\text{CH}^k(X,m;{\BQ})/ G)$ is the smallest integer
$r\ge 0$ such that there exists a closed subvariety
$i:Y \hookrightarrow X$ of [pure] codimension $k-r-m$ satisfying
$\text{CH}^k(X,m)_{\BQ} = G + i_{\ast}\text{CH}^{r+m}(Y,m;{\BQ})$.


\medskip

(ii) $\text{Level}(G)$ is the smallest integer $r\ge 0$ such that
there exists a closed subvariety
$i:Y \hookrightarrow X$ of [pure] codimension $k-r-m$ satisfying
$G \subset i_{\ast}\text{CH}^{r+m}(Y,m;{\BQ})$.
\end{defn}

Let $S$ be a smooth projective variety of dimension $r$. We
refer to the diamond below, where the upper diagonal arrows
are given by Hodge-K\"unneth projections, and the lower arrows
are defined by integration along $S$ (see \cite{Lw2}).

\begin{equation}\label{E504}
\xymatrix{& \makebox[0pt]{$\text{CH}^k(S\times X,m)$}
\ar[d] &\\
&
\makebox[36pt]{$H_{\mathcal D}^{2k-m}(S\times X,{\BR}(k)) \ar[dl]_{(m\geq 1)}$} 
\ar[dr]^{(m=0)} &\\
H^{\ell -1,0}(S)\otimes H^{k-\ell,k-m}(X) 
\ar[dr]_{\cap H^{r-\ell +1,r}(S)}
&&
H^{\ell,0}(S)\otimes H^{k-\ell,k}(X)\ar[dl]^{\cap H^{r-\ell,r}(S)}\\
&\makebox[36pt]{$H^{k-\ell,k-m}(X)$}&
}
%
%
\end{equation}

\bigskip

\begin{defn}\label{DEF505}
(i) $H^{\{k,\ell,m\}}(X) = \BC$-subspace of
$H^{k-\ell,k-m}(X)$ generated by the image of $\text{CH}^k(S\times X,m)$
in $H^{k-\ell,k-m}(X)$ in the above diagram, and over all smooth
projective algebraic $S$.

\bigskip

(ii) $H_N^{k-\ell,k-m}(X) = \BC$-subspace of $H^{k-\ell,k-m}(X)$
generated by the Hodge projected image
$$
N^{k-\ell}H^{2k-\ell -m}(X,{\BQ}) \rightarrow H^{k-\ell,k-m}(X).
$$
\end{defn}

\bigskip

\begin{rem}\label{REM506}
(i) As mentioned
in \cite{Lw2}, it is always the case
that 
\[
H^{\{k,\ell,0\}}(X) \subset H^{k-\ell,k}_N(X);
\]
moreover, one can show that $H^{\{k,\ell,0\}}(X) = 
H^{k-\ell,k}_N(X)$ under the assumption of the hard
Lefschetz conjecture. 
Further, under the assumption of the General Hodge Conjecture,
one can show that
$$
\text{Gr}_N^{k-\ell}H^{2k-\ell}(X,{\BQ})
\ne 0 \Leftrightarrow H_N^{k-\ell,k}(X) \ne 0.
$$

(ii) For $m\leq 2$, one can easily show that (see \cite{Lw2})
\[
H^{\{k-m,\ell -m,0\}}(X)\subset H^{\{k,\ell,m\}}(X).
\]
\end{rem}

\bigskip

\begin{thm}\label{THM507}
(\cite{Lw2}, abridged version) 
Let $X$ be a projective algebraic manifold and assume that $m \leq 2$ Then:
$$
H^{\{k,\ell,m\}}(X) / H^{\{k-m,\ell-m,0\}}(X) \ne 0
\Rightarrow {\rm Level}({\rm CH}^k_{\rm ind}(X,m)_{\BQ})
\geq \ell -m.
$$
\end{thm}

\begin{rem}\label{REM508}
One can readily verify that in the above theorem (see  \cite{Lw2}),
$$
\ell - m \geq 1 \Rightarrow
\text{CH}^k_\text{ind}(X,m;{\BQ}) \ \text{uncountable.}
$$

Thus if
$\text{CH}^k_\text{ind}(X,m;{\BQ})$ is uncountable,
then by rigidity, it follows that
there are an uncountable number of indecomposables in the
kernel of
$$
\text{cl}_{k,m} : \text{CH}^k(X,m;{\BQ}) \rightarrow H_{\mathcal D}^{2k-m}
(X,{\BQ}(k)),
$$
provided that the Hodge conjecture holds for $X$ in the case
$m=1$, i.e. provided that $H^{k-1,k-1}(X,{\BQ})$ is
generated by algebraic cocycles.
\end{rem}

\section{\thmref{THM507}}\label{SEC006}

It is instructive to explain some of the ideas behind \thmref{THM507}.
Broadly speaking, the the relationship between [higher] Chow
groups and Hodge theory is fortified by the following
beautiful [generalized] conjectural formula of Beilinson:
$$
\text{Gr}_F^{\ell}\text{CH}^k(X,m;{\BQ}) \simeq 
\text{Ext}^{\ell}_{\mathcal{MM}}
({\bf 1},h^{2k-\ell -m}(X)(k)),
$$
where ${\mathcal{MM}}$ is some conjectural category of mixed
motives, ${\bf 1}$ is the trivial motive, and 
\(
\text{Gr}_F^{\ell}\text{CH}^k(X,m;{\BQ})
\)
is the graded piece
of a conjectured \lq\lq Bloch-Beilinson\rq\rq\ filtration 
on $\text{CH}^k(X,m;{\BQ})$.
One way to try to realize this is to construct a duality
pairing between cohomology and higher Chow groups.  
In the case $m=0$, this is Salberger's
duality, which was exploited in \cite{Sa}, and defined as follows:
One first views $X,\ S$ as defined over an algebraically closed
field $\overline{k}$ of finite
transendence degree over $\BQ$.  Let $\eta$ be the generic point of $S$ and
$L=\overline{k}(\eta)$.  
Recall 
$$
\text{CH}^k(X_L)= \lim_{{\buildrel \rightarrow\over
{U\subset S/_{\overline{k}}}}}\text{CH}^k(U\times X).
$$
We define the top row arrow by imposing commutivity below:
\[
\begin{matrix}
\text{CH}^k(U\times X)&\times &H^{2n-2k+\ell}(X)
&- - - \longrightarrow&H^{\ell}(U)&\\
\\
\text{cl}_k\downarrow & &\downarrow Pr_2^{\ast}& &\quad\quad\uparrow\int _X(-)&\\
\\
H^{2k}(U\times X)&\times &H^{2n-2k+\ell}(U\times X)&\longrightarrow
&H^{2n+\ell}(U\times X),&\\
\\
\end{matrix}
\]
where $\int _X(-)$ is defined by integration along the fibers of
$Pr_1 : U\times X \rightarrow U$.
Taking limits, we arrive at
$$
\text{CH}^k(X_L)\times H^{2n-2k+\ell}(X)\ 
{\buildrel {\langle\ ,\ \rangle}\over\longrightarrow}\ 
H^{\ell}(\BC(S))
$$
But $\text{CH}^k(S\times X)\rightarrow \text{CH}^k(X_L)$ is surjective. 
Thus the image of $\langle\ ,\ \rangle$
lies in 
$$\text{Image}(H^{\ell}(S,{\BQ})\rightarrow H^{\ell}(\BC(S)))
\simeq \frac{H^{\ell}(S,{\BQ})}{N^1H^{\ell}(S,{\BQ})}$$

Thus: 
$$
\langle \ ,\ \rangle \ : \ \text{CH}^k(X_L)\otimes 
H^{2n-2k+\ell}(X,{\BQ})\rightarrow 
\frac{H^{\ell}(S,{\BQ})}{N^1H^{\ell}(X,{\BQ})}.
$$

The easiest way 
to explain the \thmref{THM507} is to relate it to a certain pairing
between cohomology and higher Chow groups, generalizing the
Salberger duality pairing above. A generalization of that
pairing for $\text{CH}^k(X,m;{\BQ})$ appears \cite{Lw3}. We present a 
simplified version, sufficient for our needs here.
We refer to the notation in diamond diagram \eqref{E504} above,
with $m\geq 1$.
Then there is a pairing defined in the obvious way:

\begin{equation}\label{E600}
\langle\ , \ \rangle \ : \ \text{CH}^k(S\times X,m;
{\BQ})\otimes H^{n-k+\ell,n-k+m}(X)
\ \rightarrow \ H^{\ell -1,0}(S).
\end{equation}

The trick is to relate \eqref{E600} to the level of
$\text{CH}^k(X,m;{\BQ})$.
In fact, if we view $X$ and $S$ as defined over an algebraically
closed field $\overline{k}$ of finite transcendence degree
over ${\BQ}$, and choose an embedding
$L := \overline{k}(S)\hookrightarrow \BC$, and consider the
[known injective] pullback
$\text{CH}^k(X_L,m;{\BQ})\hookrightarrow \text{CH}^k(X = X/_\BC,m;{\BQ})$,
one can argue that the pairing
in \eqref{E600} is zero if 
\[
\text{Level}(\text{CH}^k(X,m;{\BQ})) < \ell -m.
\]

\bigskip

Now let $S$ and $w \in \text{CH}^k(S\times X,m)$ be
given such that the corresponding subspace in 
$H^{\{k,\ell,m\}}(X) / H^{\{k-m,\ell-m,0\}}(X) \ne 0$.
If
\[
\text{Level}(\text{CH}^k_\text{ind}(X,m;{\BQ}))
< \ell -m,
\]
then one can argue that 
$$
\langle w,-\rangle \ : \ \{H^{\{k-m,\ell -m,0\}}(X)\}^{\perp}\ \rightarrow \
H^{\ell -1,0}(S),
$$
is zero, where
\[
\begin{split}
&\quad\{H^{\{k-m,\ell -m,0\}}(X)\}^{\perp} \subset H^{n-k+\ell,n-k+m}(X)\\
&= \{ v \ | \ \langle 
v,H^{\{k-m,\ell -m,0\}}(X) \rangle \ = 0\}.
\end{split}
\]
But by Serre duality, this in turn violates the assumption that
$$H^{\{k,\ell,m\}}(X) / H^{\{k-m,\ell-m,0\}}(X) \ne 0,$$
a contradiction.


\section{Basic strategy for Hodge-$\D$ on K3}\label{SEC007}

In the next three sections, we will prove the first part of
\thmref{THM100}, i.e., Hodge-$\D$ conjecture for a general K3
surface. As mentioned at the beginning, our basic strategy is to
degenerate a general $K3$ surface to a $K3$ surface with maximum Picard
number and study the degeneration of the higher Chow cycles 
given by \eqref{E005}.

Let $X/\Delta$ be a family of $K3$ surfaces over disk $\Delta$, where
the central fiber $X_0$ has Picard number $20$. Let
$F_1, F_2, ..., F_{20}$ be the generators of $\Pic(X_0) = \BZ^{20}$.

On every fiber $X_t$ for $t\ne 0$, we construct a higher Chow cycle
$\varepsilon_t$ in the way of \eqref{E005}. We will show
that for each $\alpha$, there are some good choices of $\varepsilon_t$
such that
the limit $\lim_{t\to 0}\cl_{2,1}(\varepsilon_t) = c_1(F_\alpha)$.

Here we have to say something about taking limit of higher Chow
cycles and regulators.
Given a family of curves $\D\subset X$, even if $\D_t$ is
reduced and irreducible, the limit $\lim_{t\to 0} \D_t$ could very
well be reducible and nonreduced. So instead of working with $\D$, we
prefer to work with its stable reduction. It naturally leads to the
following definition.

Recall that a map $\phi: C\to S$ is called stable if
\begin{enumerate}
\item $C$ is of normal crossing, i.e., only has nodes as
singularities;
\item every contractible component of $C$ under $\phi$ meets the rest
of $C$ at no less than three distinct points.
\end{enumerate}
We call $\phi$ prestable if we drop the second condition.

For a surface $S$, we define $\wt{\CH}^2(S, 1)$, 
called the higher Chow group of prestable
maps to $S$, which consists of cycles of the form
\begin{equation}\label{E007}
\sum_i (f_i, \phi_i : C_i\to S)
\end{equation}
where $\phi_i: C_i\to S$ is a prestable map from a curve $C_i$ to S,
$f_i$
is a rational function on $C_i$ such that $f_i\not\equiv 0$ when
restricted to every irreducible component of $C_i$ and
\begin{equation}\label{E008}
\sum_i \sum_{M\subset C_i} \mathop{\text{div}}
\left((\phi_i)_* \left.f_i\right|_M\right) = 0,
\end{equation}
where $M$ runs over all irreducible components of $C_i$ with
$(\phi_i)_* M \ne 0$ and
$(\phi_i)_* f_i$ is defined as follows. Let $N = \phi_i(M)_\red$
be the reduced image of $M$ under $\phi_i$ and
where $(\phi_i)_* f_i |_M\in \BC(N)^*$ is defined by
\begin{equation}\label{E009}
(\phi_i)_* f_i (p) = \prod_\alpha (f_i(q_\alpha))^{m_\alpha}
\end{equation}
for $p\in N$, where $\phi_i^* p = \sum_\alpha m_\alpha q_\alpha$ with
$m_\alpha \in \BZ$ and $q_\alpha\in M$. Note that $(\phi_i)_* f_i$ is
nothing but the norm of $f_i$ under the field extension
$\BC(M)/\BC(N)$. Of course, if either $M$ or
$N$ is singular, we pass the definition to its normalization.

We have the regulator map
$\wt{\cl}_{2,1}: \wt{\CH}^2(S, 1)\to H^{1,1}(S, \BR)^\vee$ defined by
\begin{equation}\label{E010}
\wt{\cl}_{2,1}(\wt{\varepsilon})(\omega) = \frac{1}{2\pi \sqrt{-1}}
\sum_i \sum_{M\subset C_i}\int_{M} \phi_i^* \omega\log |f_i|
\end{equation}
for $\wt{\varepsilon} = \sum_i (f_i, \phi_i : C_i\to S)$ and
$\omega\in H^{1,1}(S, \BR)$, where $M$ runs over all irreducible
components of $C_i$.

There is a natural projection
$\varphi: \wt{\CH}^2(S, 1) \to \CH^2(S, 1)$ given by
\begin{equation}\label{E011}
\varphi\left(\sum_i (f_i, \phi_i : C_i\to S) \right) =
\sum_i \sum_{M\subset C_i} ((\phi_i)_* f_i, \phi_i(M)_\red)
\end{equation}
where $M$ runs over all irreducible components of $C_i$
with $(\phi_i)_* M \ne 0$.

It is not hard to check that
\begin{equation}\label{E012}
\wt{\cl}_{2,1} = \cl_{2,1} \circ \mathop{\varphi}
\end{equation}
and it is obvious that $\varphi$ is onto. So it suffices to prove
that $\wt{\cl}_{2,1}\tensor \BR$ is surjective for a general K3
surfaces.

Let $X/\Delta$ be a smooth family of projective surfaces. It is convenient
to define the relative higher Chow group $\wt{\CH}^2(X/\Delta, 1)$ of
prestable maps,
which consists of cycles in the form
\begin{equation}\label{E014}
\wt{\varepsilon} = \sum_i (f_i, \phi_i: Y_i\to X)
\end{equation}
where $\phi_i: Y_i\to X$ is a flat family of prestable maps with the
diagram
\begin{equation}\label{E015}
\xymatrix{ Y_i \ar[r]^{\phi_i} \ar[dr] & X\ar[d]\\
& \Delta}
\end{equation}
$f_i\in \BC(Y_i)^*$ flat over $\Delta$ and $\wt{\varepsilon}_t
\in \wt{\CH}^2(X_t, 1)$ for every $t\in\Delta$.

There is a good notion of taking limit of $\wt{\cl}_{2,1}$ over a
family of surfaces $X/\Delta$. Namely, we have the following.

\begin{prop}\label{PROP001}
Let $X/\Delta$ be a smooth family of projective surfaces and
$\wt{\varepsilon} \in \wt{\CH}^2(X/\Delta, 1)$. Then
\begin{equation}\label{E013}
\lim_{t\to 0} \wt{\cl}_{2,1}(\wt{\varepsilon}_t)(\omega_t)
= \wt{\cl}_{2,1} (\wt{\varepsilon}_0) (\omega_0)
\end{equation}
for any real $(1,1)$ form $\omega_t$ of $X_t$ that varies continuously
with respect to $t$.
\end{prop}

\section{Rational Curves on $K3$ Surfaces}\label{SEC008}

\subsection{Degeneration of $K3$ surfaces}\label{SEC002001}

Consider a $K3$ surface with Picard lattice
\begin{equation}\label{E100}
\begin{pmatrix}
-2 & 1\\
1 & 0
\end{pmatrix}.
\end{equation}
Such a surface $S$ can be realized as an elliptic fibration over $\P^1$ with
fiber $F$ and a unique section $C$, where $C^2 = -2$, $C\cdot F = 1$
and $F^2 = 0$. There are exactly $24$ nodal rational curves in the
linear series $|F|$. Such surfaces were used by J. Bryan and N.C. Leung
in the enumerative problems on $K3$ surfaces \cite{B-L}. We will call
$S$ a BL $K3$ surface as in \cite{C2}. One reason why such surfaces are
so useful in the study of curves on a general $K3$ surface is that it
lies on every component of the moduli space of algebraic K3
surfaces. That is, a general $K3$ surface of genus $g$ (so the primitive
class has self-intersection $2g - 2$) can be degenerated to a
BL $K3$ surface with the primitive class degenerated to $C + gF$. In
addition, a curve $D$ in the linear series $|C+gF|$ is very easy to
describe; it is the union of $C$ and $g$ elliptic ``tails'', i.e.,
$D = C\cup F_1\cup F_2\cup ... \cup F_g$, where $F_i\in |F|$.
Moreover, if $D$ is the limit of a family of rational curves, $F_i$
must be one of the 24 rational curves; $D$ could be nonreduced and it
must be if $g > 24$.

As mentioned before, we need to degenerate a $K3$ surface to one
with maximum Picard number. Let us consider a BL $K3$ surface with
Picard number $20$. The extra Picard number comes from $-2$ curves
now appearing on the singular fibers of $S\to\P^1$. Suppose that
$S$ have $s$ singular fibers $F_1, F_2, ..., F_s$ and each $F_i$
is a closed chain of $r_i$ rational curves (see \figref{FIG001}).
The Picard number is $\sum r_i - s + 2 = 20$ and the total number
of nodes of $F_i$ is $\sum r_i = 24$. Therefore $s = 6$. We will
consider $S$ with $r_i = 4$ for $i=1,2,...,6$. That is, $S\to\P^1$
is an elliptic fiberation with $6$ singular fibers with each
singular fiber a union of $4$ rational curves (see \figref{FIG002}).

\psset{unit=0.01in,linewidth=1.2pt,dash=3pt 2pt,dotsep=2pt}


\begin{figure}[ht]
\centering
\begin{pspicture}(-170,-170)(170,0)
\psline(-150, -130)(150, -130)

\rput[lt]{0}(-150, 0){
\begin{pspicture}(0,0)(60,150)
\psline(20, 0)(20, 150) \psline(10, 120)(40, 150) \psline(30,
150)(60, 120) \psline(50, 140)(50, 90) \psline(60, 110)(30, 80)
\psline(40, 80)(10, 110)
\end{pspicture}
}

\uput[270](-130, -150){{\small $F_1$}}

\uput[270](-50, -150){{\small $F_2$}}

\uput[270](110, -150){{\small $F_s$}}

\uput[180](-150, -130){{\small $C$}}

\rput[lt]{0}(-70, 0){
\begin{pspicture}(0,0)(60,150)
\psline(20, 0)(20, 150) \psarc(0,110){30}{-70}{70}
\end{pspicture}
}

\pscircle*(20, -60){3} \pscircle*(40, -60){3} \pscircle*(60,
-60){3}

\rput[lt]{0}(90, 0){
\begin{pspicture}(0,0)(60,150)
\psline(20, 0)(20, 150) \psline(10, 150)(50, 110) \psline(50,
130)(10, 90)
\end{pspicture}
}

\end{pspicture}
\caption{A BL $K3$ surface with large Picard number}\label{FIG001}
\end{figure}
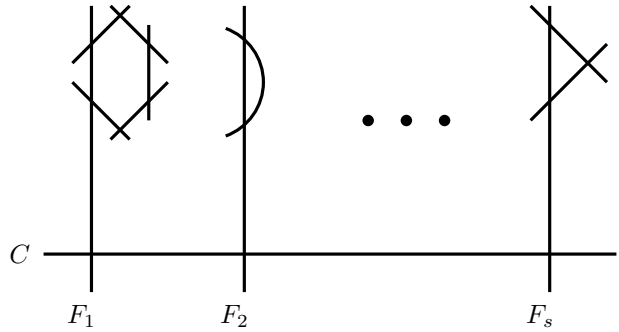


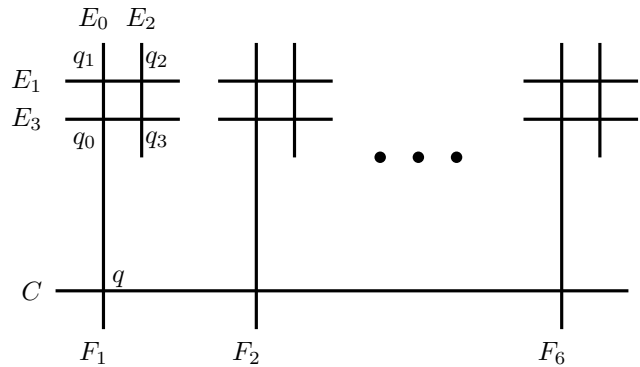
\begin{figure}[ht]
\centering
\begin{pspicture}(-170,-170)(170,0)
\psline(-150, -130)(150, -130)

\rput[lt]{0}(-150, 0){
\begin{pspicture}(0,0)(60,150)
\psline(20, 0)(20, 150) \psline(0, 130)(60, 130) \psline(40,
150)(40, 90) \psline(60, 110)(0, 110) 
\end{pspicture}
}

\uput[180](-150, -20){{\small $E_1$}}
\uput[180](-150, -40){{\small $E_3$}}
\uput[90](-130, 0){{\small $E_0$}}
\uput[90](-105, 0){{\small $E_2$}}

\uput[110](-130, -20){{\small $q_1$}}
\uput[250](-130, -40){{\small $q_0$}}
\uput[20](-110, -13){{\small $q_2$}}
\uput[340](-110, -47){{\small $q_3$}}
\uput[30](-127,-130){{\small $q$}}

\uput[270](-130, -150){{\small $F_1$}}

\uput[270](-50, -150){{\small $F_2$}}

\uput[270](110, -150){{\small $F_6$}}

\uput[180](-150, -130){{\small $C$}}

\rput[lt]{0}(-70, 0){
\begin{pspicture}(0,0)(60,150)
\psline(20, 0)(20, 150) \psline(0, 130)(60, 130) \psline(40,
150)(40, 90) \psline(60, 110)(0, 110)
\end{pspicture}
}

\pscircle*(20, -60){3} \pscircle*(40, -60){3} \pscircle*(60,
-60){3}

\rput[lt]{0}(90, 0){
\begin{pspicture}(0,0)(60,150)
\psline(20, 0)(20, 150) \psline(0, 130)(60, 130) \psline(40,
150)(40, 90) \psline(60, 110)(0, 110) 
\end{pspicture}
}

\end{pspicture}
\caption{The BL $K3$ surface we will use}\label{FIG002}
\end{figure}

Let $X/\Delta$ be a one-parameter family of $K3$ surfaces of genus $g$,
where $\Delta$ is a disk parameterized by $t$. The central fiber
$X_0 = S$ is a BL $K3$ surface described as above, given
in fig. 2.
We want to construct higher Chow cycles on $S$ using rational
curves. However, a limiting rational curve (i.e. a curve on $X_0$
which is the limit of rational curves on the general fibers $X_t$)
in $|C + gF|$ has the form $C + \sum_{i =1}^6 m_i F_i$. So any two
curves $D_1$ and $D_2$ in this form meet improperly; indeed, $D_1$
and $D_2$ do not meet properly anywhere on $S$. Thus we cannot
construct a higher Chow cycle from $D_1$ and $D_2$ in the way of
\eqref{E005}. To overcome this obstacle, we use a construction of
\cite{C2}: we blow up $X$ along $F_i$.

Let $\wt{X}$ be the blowup of $X$ along $F_1, F_2, ..., F_6$. Let us
study the behavior of the blowup along $E = F_i$.

\subsection{The blowup of $X$ along $E$}\label{SEC002002}

Let $N_{E/X}$ be the normal bundle of
$E$ in $X$. Here normal bundles are defined as the duals of
corresponding conormal bundles, as opposed to the cokernels of the
maps on the tangent spaces. We have the following exact sequence
\begin{equation}\label{E102}
\xymatrix{
0 \ar[r] & N_{E/S} \ar[r] \ar@{=}[d] & N_{E/X} \ar[r]
& \left.N_{S/X}\right|_E \ar[r] \ar@{=}[d] & 0\\
         & \CO_E & & \CO_E&
}
\end{equation}
Notice that $\Ext(\CO_E, \CO_E) = H^1(\CO_E) = \BC$ and hence the
above sequence might not split. Actually, this is always the case as
long as the family $X$ is general enough. We will sketch the argument
for this fact below. For details, please see \cite{C2}.

The long exact sequence associated to \eqref{E102} is
\begin{equation}\label{E103}
0\xrightarrow{} H^0(N_{E/S}) \xrightarrow{} H^0(N_{E/X})
\xrightarrow{} H^0(\left.N_{S/X}\right|_E) \xrightarrow{}
H^1(N_{E/S})
\end{equation}
Obviously, \eqref{E102} splits if and only if the last map
\begin{equation}\label{E104}
H^0(\left.N_{S/X}\right|_E) \xrightarrow{}
H^1(N_{E/S})
\end{equation}
is zero. We have a natural identification
$H^0(\left.N_{S/X}\right|_E)$ with $T_{\Delta, 0}$, the tangent space
of $\Delta$ at the origin. It is easy to see that \eqref{E104}
actually factors through the Kodaira-Spencer map associated to $X$,
i.e., we have
\begin{equation}\label{E105}
H^0(\left.N_{S/X}\right|_E) \isom T_{\Delta, 0}
\xrightarrow{\ks} H^1(T_S)
\xrightarrow{} H^1(T_S|_E) \isom H^1(N_{E/S}) = \BC
\end{equation}
where $\ks$ is the Kodaira-Spencer map. Here we claim that
the last map $H^1(T_S|_E)\to H^1(N_{E/S})$ is an isomorphism by the
following argument. First, we denote by $\Omega_{V}$ 
the cotangent sheaf of a variety $V$.

By the standard exact sequence
\begin{equation}\label{E800}
0\xrightarrow{} N_{E/S}^\vee \xrightarrow{} \Omega_S|_E
\xrightarrow{} \Omega_E \xrightarrow{} 0,
\end{equation}
we have the exact sequence
\begin{equation}\label{E801}
H^0(N_{E/S}) \xrightarrow{} \Ext(\Omega_E, \CO_E) \xrightarrow{}
H^1(T_S|_E)\xrightarrow{} H^1(N_{E/S}) \xrightarrow{} 0.
\end{equation}
Notice that $H^0(N_{E/S}) = \BC$ classifies the embedded deformations
of $E\subset S$ and $\Ext(\Omega_E, \CO_E) = \BC$
classifies the versal deformations of $E$. To show that $H^0(N_{E/S})$
maps nontrivially to $\Ext(\Omega_E, \CO_E)$, it suffices to show that
as $E$ varies in the pencil $|\CO_S(E)|$, the corresponding
Kodaira-Spencer map to the tangent space of
the versal deformation space of $E$ at the origin is
nontrivial, or equivalently, the map to the versal deformation space
of \(E\) is unramified over the origin.
To see this has to be true, we only need to localize the
problem at a node $p$ of $E$: if the map to the versal deformation space
is ramified over the origin, then \(S\) is locally given by \(xy =
t^\alpha\) at \(p\) for some \(\alpha > 1\); however, this is
impossible since \(S\) is smooth at \(p\). This justifies that
$H^1(T_S|_E)\to H^1(N_{E/S})$ is an isomorphism.

Since the deformation of $X_0$ in $X$ preserves the primitive class
$C + gF$, the Kodaira-Spencer class $\ks(\partial/\partial t)$ lies in
the subspace of $H^1(T_S)$ given by
\begin{equation}\label{E106}
V = \{ v\in H^1(T_S): \langle v, c_1(C + gF)\rangle  = 0\},
\end{equation}
the pairing $\langle \cdot, \cdot\rangle $ is given by Serre duality
$H^1(T_S) \times H^1(\Omega_S)\to\BC$. 
If $X$ is chosen general,
$\ks(\partial/\partial t)$ is general in $V$.

We claim that the kernel of the map
$H^1(T_S)\to H^1(T_S|_E)$ is precisely the subspace
\begin{equation}\label{E107}
W = \{ w\in H^1(T_S): \langle w, c_1(F)\rangle  = 0\}.
\end{equation}
We have the exact sequence
\begin{equation}\label{E802}
H^1(T_S(-E)) \xrightarrow{f} H^1(T_S) \xrightarrow{} H^1(T_S|_E) \isom
H^1(N_{E/S}) = \BC.
\end{equation}
So the kernel of $H^1(T_S)\to H^1(T_S|_E)$ is the image $\im f$ of $f$.
We claim that
\begin{equation}\label{E803}
\im f = W.
\end{equation}
By Kodaira-Serre duality, we have the following commutative diagram:
\begin{equation}\label{E804}
\xymatrix{H^1(T_S(-E)) \ar[r]^\sim \ar[d]^f &
H^1(\Omega_S(E))^\vee \ar[r]^\sim \ar[d] &
H^1(\Omega_S(-E)) \ar[d]^g\\
H^1(T_S)  \ar[r]^\sim & H^1(\Omega_S)^\vee \ar[r]^\sim & H^1(\Omega_S).}
\end{equation}
So we may identify the map \(f\)
with \(g: H^1(\Omega_S(-E)) \to H^1(\Omega_S)\), which is the same
as
\begin{equation}\label{E806}
g: H^{1,1}(\CO_S(-E)) \to H^{1,1}(\CO_S)
\end{equation}
on the Dolbeault cohomologies.
For any \(\psi\in H^{1, 1}(\CO_S(-E))\), we have
\begin{equation}\label{E807}
\int_S g(\psi)\wedge c_1(E) = \int_E g(\psi) = 0.
\end{equation}
So $\im f \subset W$. Since $W$ has codimension one in $H^1(T_S)$ and
\eqref{E802} is exact, we necessarily have \eqref{E803}. This
justifies our claim.

Obviously, $V\not\subset W$ and hence $\ks(\partial/\partial t)$ maps
nontrivially to $H^1(T_S|_E)$. Therefore, the map \eqref{E104} is nonzero
and \eqref{E103} does not split for $X$ general. Actually, from the
above argument we see that \eqref{E104} splits if and only if
$\ks(\partial/\partial t)\in V\cap W$, which happens if the
deformation of $X_0$ in $X$ preserves both $C$ and $F$, i.e., the
general fibers $X_t$ are also BL $K3$ surfaces.

Let $R\isom \P N_{E/X}\subset \wt{X}_0$ be the exceptional divisor
over $E$, where $\wt{X}_0$ is the central fiber of the blowup
$\wt{X}$, i.e., the total transform of $X_0$ under the blowup.
We continue to use $S$ to denote the proper transform of $S$
under $\wt{X}\to X$. The two surfaces $S$
and $R$ intersect transversely along a curve, which maps
isomorphically to $E$. Again, we continue to use $E$ for this curve.
Obviously, $E = S\cap R$ corresponds to an nonzero section in
$H^0(N_{E/X}^\vee)$. Since \eqref{E102} does not split,
$H^0(N_{E/X}^\vee) = \BC$ and
$E$ is the only section in the linear
series $\P H^0(\CO_R(E))$;
as we will see, this is the key fact which makes the geometry of $R$
interesting.

If \eqref{E102} were to split, then $R \isom \P N_{E/X}$
is simply $E\times\P^1$,
which is the {\it trivial\/} ruled surface over $E$; in our case, it
does not split so we call $R$ the {\it twisted\/} ruled surface over $E$.

Another important fact about $\wt{X}$ is that it is singular and has
exactly four rational double points over the four nodes of $E$. Let
$E = E_0\cup E_1 \cup E_2 \cup E_3$ and $q_i = E_{i-1}\cap E_i$ for
$i\in \BZ$, where we let $q_i = q_{i+4}$ and $E_i = E_{i+4}$ for
convenience (see \figref{FIG002}).
Then for each $i$, there
is a rational double point $r_i$ of $\wt{X}$ lying on the fiber of
$R\to E$ over $q_i$.

\subsection{Construction of the twisted ruled surface
$R$}\label{SEC002003}

Let $p: R\to E$ be the projection and $R_i = p^{-1} E_i$. Obviously,
$R_i\isom \F_0 := \P^1\times \P^1$.
Let $Q_i = p^{-1}(q_i)$ be the fiber
over $q_i$. An alternative way to construct $R$ is to glue four copies
of $\F_0$ along $Q_i$. The question is how to glue.

Let $\phi_{i, i+1}: Q_i\to Q_{i+1}$ be the map sending a point
$x\in Q_i$ to $y = \phi_{i, i+1}(x)\in Q_{i+1}$ such that both $x$ and
$y$ lie on a curve in the linear series $|E_i|$.
If $y = \phi_{i, i+1}(x)$, we use the notation $\overline{xy}$ to
denote the curve in $|E_i|$ passing through $x$ and $y$.
Without causing any confusion, let us abbreviate
$\phi_{i, i+1}$ to $\phi$ and
let $\phi^{-1}$ be the inverse of $\phi$.

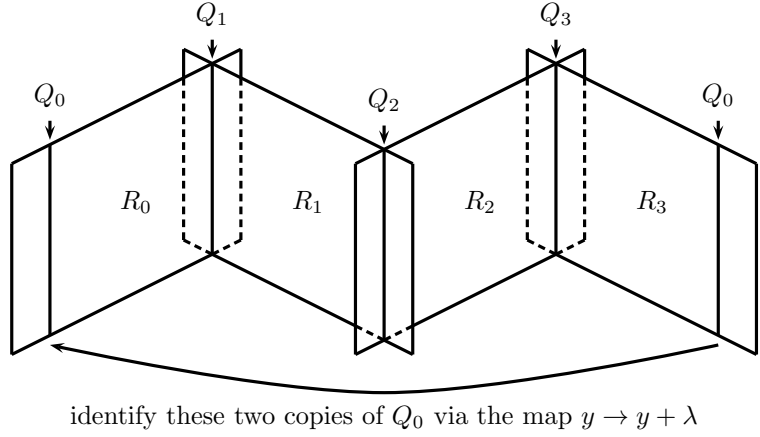
\begin{figure}[ht]
\centering
\begin{pspicture}(0, 0)(390, 215)


\psline(0, 40)(105, 92.5)
\psline[linestyle=dashed](105, 92.5)(120, 100)


\psline[linestyle=dashed](90, 100)(105, 92.5)
\psline(105, 92.5)(180, 55)
\psline[linestyle=dashed](180, 55)(195, 47.5)
\psline(195, 47.5)(210, 40)


\psline(180, 40)(195, 47.5)
\psline[linestyle=dashed](195, 47.5)(210, 55)
\psline(210, 55)(285, 92.5)
\psline[linestyle=dashed](285, 92.5)(300, 100)

\psline[linestyle=dashed](270, 100)(285, 92.5)
\psline(285, 92.5)(390, 40)

\psline(0, 140)(120, 200)
\psline(90, 200)(210, 140)
\psline(180, 140)(300, 200)
\psline(270, 200)(390, 140)

\psline(285, 92.5)(285, 192.5)
\psline{->}(285, 205)(285, 195)
\uput[90](285, 205){{\small $Q_3$}}

\psline(0, 40)(0, 140)

\psline(20, 50)(20, 150)
\psline{->}(20, 162.5)(20, 152.5)
\uput[90](20, 162.5){{\small $Q_0$}}

\uput[0](50, 120){{\small $R_0$}}
\uput[0](140,120){{\small $R_1$}}
\uput[0](230,120){{\small $R_2$}}
\uput[0](320,120){{\small $R_3$}}


\psline[linestyle=dashed](120, 100)(120, 185)
\psline(120, 185)(120, 200)

\psline(105, 92.5)(105, 192.5)

\psline{->}(105, 205)(105, 195)

\uput[90](105, 205){{\small $Q_1$}}


\psline[linestyle=dashed](90, 100)(90, 185)
\psline(90, 185)(90, 200)

\psline(210, 40)(210, 140)

\psline(195, 47.5)(195, 147.5)
\psline{->}(195, 160)(195, 150)
\uput[90](195, 160){{\small $Q_2$}}

\psline(180, 40)(180, 140)

\psline[linestyle=dashed](300, 100)(300, 185)
\psline(300, 185)(300, 200)

\psline[linestyle=dashed](270, 100)(270, 185)
\psline(270, 185)(270, 200)

\psline(390, 40)(390, 140)
\psline(370, 50)(370, 150)
\psline{->}(370, 162.5)(370, 152.5)
\uput[90](370, 162.5){{\small $Q_0$}}

\pscurve{->}(370, 45)(195, 20)(20, 45)
\uput[270](195, 20){{\small identify these two copies of $Q_0$ via the map
$y\to y+\lambda$}}

\end{pspicture}
\caption{Glue $R_i$ to obtain $R$}\label{FIG003}
\end{figure}

Consider $\phi^4: Q_0\to Q_0$. For a point $x_0\in Q_0$, let
$x_k = \phi^k(x_0)$. If $x_4 = x_0$, then
$\sum_{i=1}^4 \overline{x_{i-1} x_i}\in |E|$; however, $E$ is the only
member in the linear series $|E|$. Therefore, $x_4 = x_0$ if and only
if $x_0 = q_0$. Consequently, $\phi^4 \in \Aut(Q_0) \isom PGL(2)$ has
exactly one fixed point. If we represent $\phi^4$ by a $2\times 2$
matrix, then this matrix has exactly one eigenvector and is hence
equivalent to
\begin{equation}\label{E108}
\begin{pmatrix}
1 & \lambda\\
0 & 1
\end{pmatrix}
\end{equation}
for some $\lambda \ne 0$. More explicitly, if we pick the coordinates
on $Q_0\isom \P^1$ such that $q_0 = \infty$ becomes the point at
infinity, then $\phi^4$ is given by
\begin{equation}\label{E109}
\phi^4(y) = y + \lambda.
\end{equation}
The figure below (\figref{FIG003})
shows how to glue four copies of $\F_0$ to
obtain $R$.


Let $r_i^{(k)} = \phi^k (r_i)$ for $k\in \BZ$.
Consider the six points
\begin{equation}\label{E101}
q_0, r_0, r_0^{(4)}, r_1^{(3)}, r_2^{(2)}, r_3^{(1)}\in Q_0\isom \P^1.
\end{equation}
One can think of
$(q_0, r_0, r_0^{(4)}, r_1^{(3)}, r_2^{(2)}, r_3^{(1)})\in \m_{0, 6}$,
where $\m_{0, n}$ is the moduli space of $n$ ordered
points on $\P^1$. We claim that the moduli of these six points is
determined by the Kodaira-Spencer class $\ks(\partial/\partial t)$
of $X$ and if $\ks(\partial/\partial t)$ is general in $V$ given by
\eqref{E106}, then the corresponding point in $\m_{0, 6}$ is general.
Namely, we have the following lemma, which we will prove later.

\begin{lem}\label{LEM001}
There is a well-defined rational map $\P V\to \m_{0, 6}$ which sends
\begin{equation}\label{E901}
v \to (q_0, r_0, r_0^{(4)}, r_1^{(3)}, r_2^{(2)}, r_3^{(1)})
\end{equation}
where $q_0, r_i, \phi$ are obtained
from a family $X$ with Kodaira-Spencer class $v$.
Furthermore $\P V$ dominates $\m_{0, 6}$ under this map.
\end{lem}

By the above lemma, we see that $r_i^{(k)}\ne r_j^{(l)}$ for all
$k,l\in \BZ$ provided $i\ne j$.

\subsection{Construction of Limit Rational Curves on $\wt{X}_0$}
\label{SEC002004}

Let $q = C\cap E$, where we recall that $C$ is the unique section
$S\to \P^1$ with $C^2 = -2$. Without loss of generality, we assume
that $q\in E_0$.

Let $\Gamma$ be a curve on $R$ satisfying the following conditions.

\begin{enumerate}
\item $\Gamma$ is reduced and $E_i, Q_i\not\subset \Gamma$ for all
$E_i$ and $Q_i$.
\item Let $\Gamma_i = \Gamma\cap R_i$. Then
$\Gamma_i \in |m E_i|$ for $R_i\ne R_0$; $\Gamma_0 = N\cup G$, where
$N$ is an irreducible curve in $|Q_0 + \mu E_0|$,
$q\in N$ and $G\in |(m - \mu) E_0|$, and where $m$, $\mu$ are
non-negative integers.
\item $\Gamma$ consists of $N$ and chains of curves in the form
\begin{equation}\label{E110}
\sum_{j = k}^{l-1} \overline{r_i^{(j)} r_i^{(j+1)}}
\end{equation}
where $k\le 0$, $l\ge 0$, $r_i^{(k)}\in Q_1$, $r_i^{(l)} \in Q_0$
(i.e. $i \equiv 1 - k\equiv -l (\mod 4)$) and $N$ passes through
$r_i^{(k)}$ and $r_i^{(l)}$.
\item $N$ meets $Q_0$ and $Q_1$ transversely everywhere and if $N$
meets $Q_0$ or $Q_1$ at some point $u$,
then $\Gamma$ contains a chain of curves in the form of \eqref{E110}
with $u = r_i^{(k)}$ or $r_i^{(l)}$.
\end{enumerate}

For examples of such $\Gamma$, please see the figures in the next
section.

We claim that such $\Gamma$ is part of a limiting rational
curve. More precisely, we have the following.

\begin{prop}\label{PROP005}
Let $D$ be a curve on $\wt{X}_0$ with $\nu_* D \in |C + gF|$, where
$\nu$ is the blowup map $\wt{X}\to X$. Assume that for each $E = F_s$
($s = 1,2,...,6$), the corresponding $\Gamma = D\cap R$ satisfies all
the conditions (1)-(4) listed above.
Then $D$ can be deformed to rational curves
on the general fibers $\wt{X}_t$.
More precisely, there exists a family of stable rational maps
$\wt{\pi}: Y\to \wt{X}$ with the diagram
\begin{equation}\label{E113}
\xymatrix{Y \ar[r]^{\wt{\pi}} \ar[dr]^{\pi} \ar[ddr] & \wt{X}
\ar[d]\\ & X \ar[d]\\
& \Delta}
\end{equation}
such that $D = \wt{\pi}_* Y_0$.
\end{prop}

The proof of the above proposition is not hard since $D$ is
reduced and has only nodes as singularities. We just have to figure
out when $D$ deforms, which nodes of $D$ remain as nodes and which are
smoothed out. Using the deformational argument in \cite{C1} and
\cite{C2}, we see that $N\cap G$ and $r_i\in \Gamma$ remain as nodes
and the rest are smoothed out. As a consequence, we can describe $Y_0$
as follows.

\begin{prop}\label{PROP004}
Let $Y, D, \Gamma$ be given as in \propref{PROP005} and
let $\wt{\Gamma}\subset Y_0$ be the pre-image of $\Gamma$. Then
$\wt{\Gamma}$ meets the rest of $Y_0$ at a point over $q$. And
\begin{enumerate}
\item $\wt{\Gamma}$ consists of $\wt{N}$, which dominates $N$,
and chains of rational curves attached to $\wt{N}$;
\item each chain of curves \eqref{E110} on $\Gamma$ breaks
into (at most) two chains of curves on $\wt{\Gamma}$: one dominates
\begin{equation}\label{E111}
\sum_{j = k}^{-1} \overline{r_i^{(j)} r_i^{(j+1)}}
\end{equation}
the other dominates
\begin{equation}\label{E112}
\sum_{j = 0}^{l-1} \overline{r_i^{(j)} r_i^{(j+1)}}
\end{equation}
and they meet $\wt{N}$ at points over $r_i^{(k)}$ or $r_i^{(l)}$.
\end{enumerate}
\end{prop}

We certainly did not give all possible limiting rational curves in
\propref{PROP005}.
The curve $D$ described there only represents a small fraction of all
possible degenerations of rational curves on the general fiber. Using
the argument in \cite{C2}, one can classify all limiting rational
curves. But there is no need for that here. We only need those $D$'s
decribed above.

\section{Proof of Hodge-$\D$-conjecture for General $K3$ Surfaces}
\label{SEC009}

\subsection{Construction of higher Chow cycles}\label{SEC003001}

Consider a pair of curves $\Gamma, \Sigma\subset R$,
satisfying (1)-(4) in \ssecref{SEC002004}:
\begin{equation}\label{E300}
\Gamma = N_\Gamma + \sum_{j = -3}^{-1} \overline{r_0^{(j)}r_0^{(j+1)}}
\end{equation}
where $N_\Gamma \subset R_0$ is the unique curve in
$|E_0 + Q_0|$ passing through $r_0^{(-3)}, r_0, q$;
\begin{equation}\label{E301}
\Sigma =  N_\Sigma + \sum_{j = 0}^2 \overline{r_1^{(j)} r_1^{(j+1)}}
+ \sum_{j = -1}^1 \overline{r_2^{(j)} r_2^{(j+1)}}
\end{equation}
where $N_\Sigma\subset R_0$ is the unique curve in
$|E_0 + 2 Q_0|$ passing through the five points
$r_1, r_1^{(3)}, r_2^{(-1)}, r_2^{(2)}, q$
(see \figref{FIG004} and \ref{FIG005}).

Obviously, $\Gamma$ and $\Sigma$ meet at three points
$\{u, v, q\} = N_\Gamma\cap N_\Sigma$.

By \propref{PROP005}, there exists families of stable rational maps
$\wt{\pi}_Y: Y\to \wt{X}$ and $\wt{\pi}_Z: Z\to \wt{X}$ with the
diagram
\begin{equation}\label{E900}
\xymatrix{Y \ar[r]^{\wt{\pi}_Y} \ar[dr]^{\pi_Y} \ar[ddr] & \wt{X}
\ar[d] & Z \ar[l]_{\wt{\pi}_Z} \ar[dl]_{\pi_Z} \ar[ddl]\\
& X \ar[d]\\
& \Delta}
\end{equation}
such that
$(\pi_Y)_* Y_0, (\pi_Z)_* Z_0\in |C + gF|$,
$\Gamma = R\cap (\wt{\pi}_Y)_* Y_0$ and
$\Sigma = R\cap (\wt{\pi}_Z)_* Z_0$.

\figref{FIG004} shows $\Gamma$ and its pre-image on $Y_0$
and \figref{FIG005} shows $\Sigma$ and its pre-image on $Z_0$.


\begin{figure}[ht]
\centering
\begin{pspicture}(0, 0)(400, 150)
\pspolygon(0, 20)(0, 120)(60, 120)(60,20)
\psline(60, 120)(100, 140)(100, 40)(60, 20)
\psline(100, 140)(140, 120)(140, 20)(100, 40)
\psline(140, 120)(180, 140)(180, 40)(140, 20)

\uput[90](30, 120){{\small $R_0$}}
\uput[90](80, 130){{\small $R_1$}}
\uput[90](120,130){{\small $R_2$}}
\uput[90](160,130){{\small $R_3$}}

\psline(60, 90)(100, 110)(140, 90)(180, 110)

\uput[0](180, 110){{\small $r_0$}}
\uput[0](20, 70){{\small $N_\Gamma$}}

\pscurve(36, 60)(40, 35)(50, 20)
\pscurve(60, 90)(40, 100)(36, 80)
\pscurve(0, 110)(10, 100)(25, 70)

\uput[270](50, 20){{\small $q$}}
\uput[180](0, 110){{\small $r_0$}}

\psline{->}(190, 90)(210, 90)

\psline(210, 20)(400, 20)
\uput[0](400,20){{\small $C$}}
\uput[270](305, 15){{\small on $Y_0$}}

\psline(230, 0)(230, 140)
\uput[90](230, 140){{\small $N_\Gamma$}}
\uput[225](230, 20){{\small $q$}}

\psline(220, 100)(260, 140)
\psline(250, 140)(290, 100)
\psline(280, 100)(320, 140)

\uput[300](240, 120){{\small $E_1$}}
\uput[60](270, 120){{\small $E_2$}}
\uput[300](300, 120){{\small $E_3$}}

\pscircle*(310, 130){2}
\psline{<-}(320, 130)(340, 130)
\uput[0](340, 130){{\small $r_0$}}

\end{pspicture}
\caption{Curve $\Gamma$}\label{FIG004}
\end{figure}


\begin{figure}[ht]
\centering
\begin{pspicture}(0, 0)(400, 160)
\pspolygon(0, 20)(0, 120)(60, 120)(60,20)
\psline(60, 120)(100, 140)(100, 40)(60, 20)
\psline(100, 140)(140, 120)(140, 20)(100, 40)
\psline(140, 120)(180, 140)(180, 40)(140, 20)

\uput[90](30, 120){{\small $R_0$}}
\uput[90](80, 130){{\small $R_1$}}
\uput[90](120,130){{\small $R_2$}}
\uput[90](160,130){{\small $R_3$}}

\psline(60, 90)(100, 110)(140, 90)(180, 110)

\uput[320](60,90){{\small $r_1$}}

\psline(60, 50)(100, 70)(140, 50)(180, 70)

\uput[30](100,70){{\small $r_2$}}

\uput[0](180, 110){{\small $r_1^{(3)}$}}
\uput[0](20, 70){{\small $N_\Sigma$}}

\pscurve(36, 60)(40, 35)(50, 20)
\pscurve(60, 90)(40, 100)(36, 80)
\pscurve(0, 110)(10, 100)(25, 70)
\psline(60, 50)(40, 60)
\pscurve(0, 70)(10, 55)(25,60)

\uput[180](0, 70){{\small $r_2^{(2)}$}}
\uput[270](50, 20){{\small $q$}}
\uput[180](0, 110){{\small $r_1^{(3)}$}}

\psline{->}(190, 90)(210, 90)

\psline(210, 20)(400, 20)
\uput[0](400,20){{\small $C$}}
\uput[270](305, 15){{\small on $Z_0$}}

\psline(230, 0)(230, 140)
\uput[90](230, 140){{\small $N_\Sigma$}}
\uput[225](230, 20){{\small $q$}}

\psline(220, 100)(260, 140)
\psline(250, 140)(290, 100)
\psline(280, 100)(320, 140)

\uput[300](240, 125){{\small $E_3$}}
\uput[60](270, 115){{\small $E_2$}}
\uput[300](300, 125){{\small $E_1$}}

\pscircle*(310, 130){2}
\psline{<-}(320, 130)(340, 130)
\uput[0](340, 130){{\small $r_1$}}

\psline(220, 60)(260, 100)
\psline(250, 100)(290, 60)

\uput[300](240, 85){{\small $E_3$}}
\uput[60](270, 75){{\small $E_2$}}

\pscircle*(280, 70){2}
\psline{<-}(290, 70)(310, 70)
\uput[0](310, 70){{\small $r_2$}}

\psline(220, 50)(270, 50)
\uput[270](245, 50){{\small $E_1$}}
\pscircle*(260, 50){2}

\pscurve{->}(315,65)(275, 40)(262, 48)
\end{pspicture}
\caption{Curve $\Sigma$}\label{FIG005}
\end{figure}
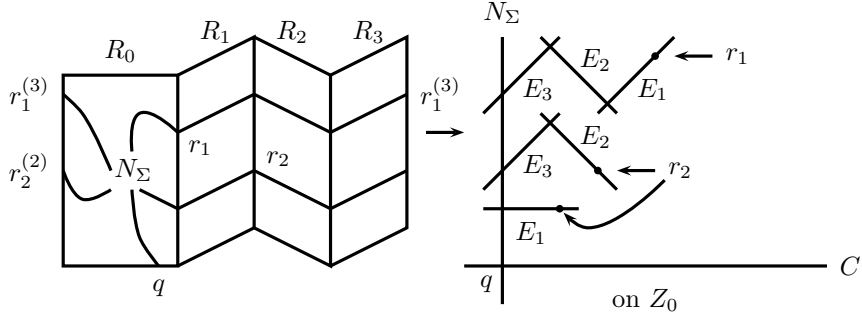

Also $\wt{\pi}_Y(Y_0)$ and $\wt{\pi}_Z(Z_0)$ meet transversely
at $u$ and $v$. Therefore, $\wt{\pi}_Y(Y_t)$ and $\wt{\pi}_Z(Z_t)$
meet at two points in the neighborhood of $u$ and $v$, respectively.
That is, there exist sections
$U_Y\subset Y$ of $Y/\Delta$ and
$U_Z \subset Z$ of $Z/\Delta$ such that
$u\in \wt{\pi}_Y(U_Y) = \wt{\pi}_Z(U_Z)$; there exist sections
$V_Y\subset Y$ of $Y/\Delta$ and
$V_Z \subset Z$ of $Z/\Delta$ such that
$v\in \wt{\pi}_Y(V_Y) = \wt{\pi}_Z(V_Z)$.

We construct a cycle
\begin{equation}\label{E302}
\wt{\varepsilon} = (f_Y, \pi_Y: Y\to X) +
(f_Z, \pi_Z: Z\to X)
\end{equation}
of $\wt{\CH}^2(X/\Delta, 1)$ with
\begin{equation}\label{E303}
(f_Y) = (U_Y) - (V_Y) \text{ and } (f_Z) = (V_Z) - (U_Z).
\end{equation}

\subsection{Computation of $\wt{\varepsilon}_0$}\label{SEC003002}

Thanks to our knowledge about $R$ obtained in \ssecref{SEC002003},
we can make
explicit calculation of $\wt{\varepsilon}_0$.

Let $R_0 = E_0\times Q_0$ be parametrized by affine coordinates
$(x, y)$ where $x$ and $y$ are the affine coordinates of $E_0$ and
$Q_0$, respectively. We choose $(x, y)$ such that
$E_0 = \{y = \infty\}$, $Q_0 = \{ x = 0 \}$, $Q_1 = \{ x = \infty \}$
and $q = (1, \infty)$.

Assume $r_1 = (\infty, y_1)$, $r_2^{(-1)} = (\infty, y_2)$
and $r_0^{(-3)} = (\infty, y_0)$. And assume that $\phi^4$ is
given by \eqref{E109}. We may choose $(x, y)$ such that
$\lambda = 1$.
Then $r_1^{(3)} = (0, y_1 + 1)$,
$r_2^{(2)}  = (0, y_2 + 1)$ and $r_0 = (0, y_0 + 1)$.

Then $N_\Gamma$ is given by
\begin{equation}\label{E304}
x (y - y_0) - (y - y_0 - 1) = 0
\end{equation}
and $N_\Sigma$ is given by
\begin{equation}\label{E305}
x (y - y_1)(y - y_2) - (y - y_1 - 1)(y - y_2 - 1) = 0.
\end{equation}
Then the $y$-coordinates of $u$ and $v$ are the two roots of the
equation
\begin{equation}\label{E306}
\frac{y - y_0 - 1}{y - y_0} = \frac{(y - y_1 - 1)(y - y_2 - 1)}{(y -
y_1)(y - y_2)}.
\end{equation}
Let $\alpha_0$ and $\beta_0$ be the roots of the above equation.
Then the restrictions of $f_Y$ and $f_Z$ to $N_\Gamma$ and $N_\Sigma$
are
\begin{equation}\label{E307}
g_0(y) = \frac{y - \alpha_0}{y - \beta_0} \text{ and }
\frac{1}{g_0(y)} = \frac{y - \beta_0}{y - \alpha_0}
\end{equation}
and the restrictions of $f_Y$ and $f_Z$ to the rest of $Y_0$ and
$Z_0$ are constants (see \figref{FIG004} and \ref{FIG005}).

So $\wt{\varepsilon}_0$ looks roughly like this
\begin{equation}\label{E808}
\begin{split}
\wt{\varepsilon}_0 &= (g_0, N_\Gamma) + \sum (\text{constant},
\text{other component of $Y_0$})\\
& + (\frac{1}{g_0}, N_\Sigma) +
\sum (\text{constant}, \text{other component of $Z_0$}),
\end{split}
\end{equation}
where the ``constants'' are the values of $g_0$ and $1/g_0$ at the
intersections between $N_\Gamma, N_\Sigma$ and the other components of
$Y_0, Z_0$. Rigorously, we should write e.g. $(g_0, N_\Gamma)$ as
$(g_0, \pi_Y: N_\Gamma\to X_0)$ but we want to keep our expression
simple.

>From \figref{FIG004}, $N_\Gamma$ meets $E_1\cup E_2\cup E_3$ at a
point over $r_0^{(-3)} = (\infty, y_0)$ and the value of $g_0$ at
$r_0^{(-3)}$ is $g(y_0)$. So the restriction of $f_Y$ to
$E_1\cup E_2\cup E_3$ is the constant $g(y_0)$.
Hence $\wt{\varepsilon}_0$ contains the cycle
\begin{equation}\label{E809}
(g(y_0), E_1 + E_2 + E_3).
\end{equation}
Rigorously, $E_i$ ($i=1,2,3$) is a component of $Y_0$ dominating
$E_i\subset X_0$. But we would not use different notations as it would
just make our answer a lot messier than necessary.
Also notice that $N_\Gamma$ meets $C$ at $q = (1, \infty)$,
where $g_0$ has value
$1$. So the restriction of $f_Y$ to the rest of $Y_0$ is $1$. That is,
these cycles are killed under the regulator map.

>From \figref{FIG005}, $N_\Sigma$ meets $E_3\cup E_2\cup E_1$ at a
point over $r_1^{(3)} = (0, y_1 + 1)$. Hence $\wt{\varepsilon}_0$
contains the cycle
\begin{equation}\label{E810}
-(g(y_1 + 1), E_3 + E_2 + E_1).
\end{equation}
Here, for convenience, we extend the group law of $\CH^2(S, 1)$
by letting $(f, D_1) + (f, D_2) = (f, D_1 + D_2)$ and
$n(f, D) = (f^n, D)$ for $n\in \BZ$. Obviously, this is compatible
with the regulator map.

Similarly, $N_\Sigma$ meets $E_3\cup E_2$ at a point over
$r_2^{(2)} = (0, y_2 + 1)$ and $N_\Sigma$ meets $E_1$ at a point over
$r_2^{(-1)} = (\infty, y_2)$. Hence $\wt{\varepsilon}_0$
contains the cycles
\begin{equation}\label{E811}
-(g(y_2+1), E_3 + E_2) - (g(y_2), E_1).
\end{equation}
Combining \eqref{E809}-\eqref{E811}, we have
\begin{equation}\label{E812}
\begin{split}
\wt{\varepsilon}_0 &= (g_0, N_\Gamma) + (g(y_0), E_1 + E_2 + E_3)\\
& -(g_0, N_\Sigma) - (g(y_1 + 1), E_3 + E_2 + E_1)\\
&-(g(y_2+1), E_3 + E_2) - (g(y_2), E_1) + (1, *)
\end{split}
\end{equation}
where by $(1, *)$ we mean that the rest terms are
of the form $(1, D)$.

Recall that $\varphi$ is the projection
$\wt{\CH}^2(S, 1)\to \CH^2(S, 1)$ defined in \secref{SEC007},
where $S = X_0$. Now we push $\wt{\varepsilon}_0$ forward to
$\varepsilon_0 = \varphi(\wt{\varepsilon}_0)$.
Note that
\begin{equation}\label{E813}
\varphi(g_0, N_\Gamma) = ((\pi_Y)_* g_0, E_0)\ \text{and}\
\varphi(g_0, N_\Sigma) = ((\pi_Z)_* g_0, E_0).
\end{equation}
Obviously, $\text{div} (\pi_Y)_* g_0 = \text{div} (\pi_Z)_* g_0$. Therefore,
$((\pi_Y)_* g_0) / ((\pi_Z)_* g_0)$ is a constant. To find this
constant, it is enough to evaluate $(\pi_Y)_* g_0$ and $(\pi_Z)_* g_0$
at some point on $E_0$, say $q$. Since $\pi_Y: N_\Gamma \to E_0$ is
one-to-one, the preimage of $q$ is itself and hence $(\pi_Y)_* g_0(q)
= g_0(q) = 1$. And $\pi_Z: N_\Sigma\to E_0$ is two-to-one, the preimage of
$q$ consists of $q$ and a point with $y$-coordinate
$(y_1 + y_2 + 1)/2$ (solve \eqref{E305} after setting
$x = 1$). Therefore,
\begin{equation}\label{E814}
\frac{(\pi_Y)_* g_0}{(\pi_Z)_* g_0} = \frac{1}{g_0((y_1 + y_2 + 1)/2)}.
\end{equation}

In conclusion, we have
\begin{equation}\label{E308}
\begin{split}
\varepsilon_0 = & -\left(g_0(\frac{y_1 + y_2 + 1}{2}),
E_0\right)\\
& +
\left(\frac{g_0(y_0)}{g_0(y_1+1)g_0(y_2 + 1)},
E_1 + E_2 + E_3\right)\\
& + \left(\frac{g_0(y_2 + 1)}{g_0(y_2)}, E_1\right) + (1, *).
\end{split}
\end{equation}

Therefore, $\cl_{2,1}(\varepsilon_0)$ is
\begin{equation}\label{E309}
\begin{split}
& -\log \left|g_0(\frac{y_1 + y_2 +
1}{2})\right| c_1(E_0)\\
& + \log\left|\frac{g_0(y_0)}{g_0(y_1+1)g_0(y_2 + 1)}\right| c_1(E_1 +
E_2 + E_3)\\
&+ \log \left|\frac{g_0(y_2 + 1)}{g_0(y_2)}\right| c_1(E_1).
\end{split}
\end{equation}

Next, we will change $\Gamma$ and $\Sigma$ to produce more classes
in $H^{1, 1}(S, \BR)^\vee$. Let $\Gamma$ and $\Sigma$ be given as
follows:
\begin{equation}\label{E310}
\Gamma = N_\Gamma + \sum_{j=0}^2 \overline{r_1^{(j)} r_1^{(j+1)}}
\end{equation}
where $N_\Gamma\subset R_0$ is the unique curve in $|E_0 + Q_0|$
passing through $r_1, r_1^{(3)}, q$ and
\begin{equation}\label{E311}
\Sigma = N_\Sigma + \sum_{j = -3}^{-1} \overline{r_0^{(j)} r_0^{(j+1)}}
+ \sum_{j = -1}^{1} \overline{r_2^{(j)} r_2^{(j+1)}}
\end{equation}
where $N_\Sigma\subset R_0$ is the unique curve in $|E_0 + 2Q_0|$
passing through the four points $r_0^{(-3)}, r_0, r_2^{(-1)}, r_2^{(2)}, q$
(compare \eqref{E310} and \eqref{E311} with \eqref{E300} and
\eqref{E301}).

The corresponding $\cl_{2,1}(\varepsilon_0)$ is
\begin{equation}\label{E312}
\begin{split}
&-\log\left|g_1(\frac{y_2 + y_0 + 1}{2})\right| c_1(E_0)\\
& + \log \left |\frac{g_1(y_1 + 1)}{g_1 (y_0) g_1(y_2+1)}\right|
c_1(E_1 +
E_2 + E_3)\\
&+ \log \left|\frac{g_1(y_2 + 1)}{g_1(y_2)}\right| c_1(E_1)
\end{split}
\end{equation}
where
\begin{equation}\label{E313}
g_1(y) = \frac{y - \alpha_1}{y - \beta_1}
\end{equation}
and $\alpha_1$ and $\beta_1$ are the two roots of
\begin{equation}\label{E314}
\frac{y - y_1 - 1}{y - y_1} = \frac{(y - y_2 - 1)(y - y_0 - 1)}{(y -
y_2)(y - y_0)}.
\end{equation}

We may produce one more class by choosing $\Gamma$ and $\Sigma$ to be
\begin{equation}\label{E315}
\Gamma = N_\Gamma + \sum_{j=-1}^1 \overline{r_2^{(j)} r_2^{(j+1)}}
\end{equation}
where $N_\Gamma\subset R_0$ is the unique curve in $|E_0 + Q_0|$
passing through $r_2^{(-1)}, r_1^{(2)}, q$ and
\begin{equation}\label{E316}
\Sigma = N_\Sigma + \sum_{j = -3}^{-1} \overline{r_0^{(j)} r_0^{(j+1)}}
+ \sum_{j = 0}^{2} \overline{r_1^{(j)} r_1^{(j+1)}}
\end{equation}
where $N_\Sigma\subset R_0$ is the unique curve in $|E_0 + 2Q_0|$
passing through the four points $r_0^{(-3)}, r_0, r_1, r_1^{(3)}, q$.

The corresponding $\cl_{2,1}(\varepsilon_0)$ is
\begin{equation}\label{E317}
\begin{split}
&-\log\left|g_2(\frac{y_0 + y_1 + 1}{2})\right| c_1(E_0)\\
& + \log \left |\frac{g_2(y_2 + 1)}{g_2 (y_0) g_2(y_1+1)}\right|
c_1(E_1 +
E_2 + E_3)\\
&- \log \left|\frac{g_2(y_2 + 1)}{g_2(y_2)}\right| c_1(E_1)
\end{split}
\end{equation}
where
\begin{equation}\label{E318}
g_2(y) = \frac{y - \alpha_2}{y - \beta_2}
\end{equation}
and $\alpha_2$ and $\beta_2$ are the two roots of
\begin{equation}\label{E319}
\frac{y - y_2 - 1}{y - y_2} = \frac{(y - y_0 - 1)(y - y_1 - 1)}{(y -
y_0)(y - y_1)}.
\end{equation}

Finally, we need to show that the $3\times 3$ matrix
formed by the coefficients of \eqref{E309}, \eqref{E312} and
\eqref{E317} is invertible. This is easy thanks to \lemref{LEM001},
which tells us that $y_0, y_1, y_2$ can be chosen arbitrarily. With
the help of a computer, one can easily find some numerical values of
$y_0, y_1, y_2$ for which the matrix is invertible. For example, we
find $y_0 = 0, y_1 = 1/8, y_2 = 1/2$ a good choice. The Maple program
we used is available upon request.

Therefore, $c_1(E_0)$, $c_1(E_1 + E_2 + E_3)$ and $c_1(E_1)$ are in
the image of $\cl_{2,1}\tensor \BR$. Change $r_2$ to $r_4$ and by the
same argument we will obtain $c_1(E_3)$ in the image of
$\cl_{2,1}\tensor \BR$. Since $E$ can be chosen to be any of the six
singular fibers $F_1, F_2, ..., F_6$, we obtain all the classes in
$\Pic(S)$ and we are done. It only remains to prove \lemref{LEM001}.

\subsection{Proof of \lemref{LEM001}}\label{SEC003003}

We have the surjective map
\begin{equation}\label{E320}
N_{E/X} \xrightarrow{} \ext^1(\Omega_E, \CO_E) \xrightarrow{} 0
\end{equation}
where $\ext^1(\Omega_E, \CO_E) = T^1(E)$
is usually called the $T^1$ sheaf of $E$, which is a sheaf supported
on the singular locus of $E$, i.e., $q_i$. It is easy to see that
$r_i\in Q_i$ is induced by the map
\begin{equation}\label{E321}
\left. N_{E/X}\right |_{q_i} \to \left. T^1(E) \right|_{q_i}.
\end{equation}
Therefore, we see that $(q_0, r_0, \phi^4 (r_0), \phi^3(r_1), \phi^2
(r_2), \phi (r_3))\in \m_{0,6}$ only depends on the Kodaira-Spencer
class of $X$ and the map \eqref{E901} is well defined.

Instead of $\m_{0,6}$, we consider the map $f: \P V \to (\P^1)^3$
sending $v\in V$ to $(\phi^3(r_1), \phi^2 (r_2), \phi (r_3))$ with
$(q_0, r_0, \phi^4 (r_0))$ fixed at $(0, 1, \infty)$. Obviously, $f$
is dominant if and only if the original map $\P V\to \m_{0,6}$ is
dominant.

Let $V_3\subset V_2\subset V_1 \subset V_0 = V$ be a filtration,
where
\begin{equation}\label{E322}
\begin{split}
V_i = \{v\in V: \langle v, c_1(E_1)\rangle  &= \langle v, c_1(E_2)\rangle  \\
&= ... = \langle v, c_1(E_i)\rangle 
= 0\}.
\end{split}
\end{equation}

We claim the following
\begin{claim}\label{CLM001}
$r_i = \phi(r_{i+1})$ if and only if $\langle v, c_1(E_i)\rangle  = 0$.
\end{claim}

If the above claim is true, then we have
\begin{equation}\label{E323}
f(\P V_3) \subsetneq f(\P V_2) \subsetneq f(\P V_1) \subsetneq
f(\P V_0)
\end{equation}
and it follows that $\dim f(\P V_0) = 3$, i.e., $f$ is dominant.
So it remains to justify our claim.

Without the loss of generality, take $E_i = E_1$.
Let $Z$ be the blowup of $X$ along $E_1$ with the exceptional divisor
$M \isom \P N_{E_1/X}$. By the same argument as in
\ssecref{SEC002002}, we can show that the exact sequence
\begin{equation}\label{E324}
0 \xrightarrow{} N_{E_1/S} \xrightarrow{} N_{E_1/X}
\xrightarrow{} \left. N_{S/X} \right |_{E_1} \xrightarrow{} 0
\end{equation}
splits if and only if $\langle v, c_1(E_1)\rangle  = 0$.
Therefore, $M\isom \F_2$ if $\langle v, c_1(E_1)\rangle  = 0$ and
$M\isom \F_0 = \P^1\times \P^1$ otherwise.

Let $\wt{Z}$ be the blowup of $Z$ along the proper transform of
$E_0\cup E_2\cup E_3$ under $Z\to X$. Then $\wt{Z}$ is actually the
resolution of $\wt{X}$
at the rational double points $r_1$ and $r_2$. We
have the diagram
\begin{equation}\label{E325}
\xymatrix{\wt{Z} \ar[r] \ar[d] & X \ar[ddl]\\
Z \ar[d]\\
X}
\end{equation}
Let $\wt{M}\subset \wt{Z}$ be the pullback of $M$ under $\wt{Z}\to
Z$. Then we have the diagram
\begin{equation}\label{E326}
\xymatrix{ \wt{M} \ar[r] \ar[d] & R_1\\
M}
\end{equation}
Both the maps $\wt{M}\to M$ and $\wt{M} \to R_1$ are the blowdowns of
two $(-1)$-curves. But the two sets of $(-1)$-curves are
distinct. They can be describe explicitly as follows.

Starting with $M$, we blow up $M$ at two points $p_1$ and $p_2$ to
obtain $\wt{M}$. Obviously, $p_1$ and $p_2$ are over $q_1$ and $q_2$,
respectively.
Let $G_1$ and $G_2$ be the two fibers of $M\to E_1$
that contain $p_1$ and $p_2$, respectively, and let $\wt{G}_1$ and
$\wt{G}_2$ be the proper transforms of $G_1$ and $G_2$ under the blow
up. It is easy to see that the map $\wt{M}\to R_1$
blows down $\wt{G}_1$ and $\wt{G}_2$ to $r_1$ and $r_2$. Under such
description, one can easily see that $\phi(r_1) = r_2$
if $M\isom \F_2$ and $\phi(r_1) \ne r_2$ if $M\isom \F_0$.
[Here $\F_{n} = \P(\CO\oplus \CO(n))$.]

\section{Proof of \thmref{THM100}}\label{SEC010}

In \cite{S-S} there is constructed a $19$-dimensional
``universal'' family $S$ of marked $K3$ surfaces,
where roughly speaking, a marked $K3$ surface is
a $K3$ surface together with the even unimodular lattice
$L := H_{2}(X,{\BZ})$, and with choice of polarization
$\xi \in L$. The special Kummer surfaces form a 
dense subset of that family. However we learned from 
the previous sections that the Hodge-${\mathcal D}$-conjecture 
holds over a real Zariski open subset of that family.
Thus it is clear that the Hodge-${\mathcal D}$-conjecture holds
for general special Kummer surfaces, as well as for general Kummer 
surfaces as well. By ``lifting'' these  results to
the corresponding Abelian surfaces, we can now prove the following.
[General will mean as in \secref{SEC002},
in the real analytic Zariski topology.]

\bigskip

\begin{thm}\label{THM000}
{\rm (i)}
The Hodge-${\mathcal D}$-conjecture
holds for surfaces of the form $E_{1}\times E_{2}$,
where $E_{1},\ E_{2}$ are general elliptic curves.

\bigskip

{\rm (ii)}
The Hodge-${\mathcal D}$-conjecture
holds for general Abelian surfaces.
\end{thm}

\bigskip

\begin{proof}
Let $A$ be an Abelian surface,
$Y = A/\pm 1$, be the quotient space after\
applying the involution $x\mapsto -x$, and $X$
the corresponding Kummer surface resulting from
the blow-up of the $16$ double points in $Y$. It
is clear that the cohomology $H^{2}(X)$, which
is of rank $22$, comes from the cohomology of
$H^{2}(A)$ (of rank $6$) together with the $16$
rational curves resulting from the aforementioned
blow-up process. Let $w\subset A\times X$ be the
correspondence induced by the above process.
Then 
$$
[w]^{\ast} : H^{2}(X,{\BR}(1))\to H^{2}(A,{\BR}(1)),
$$
is a surjective morphism of Hodge structures. However
by functoriality, there is a commutative diagram:
\[
\begin{matrix}
\text{CH}^{2}(X,1;{\BR})&{\buildrel w^{\ast}\over\to}
&\text{CH}^{2}(A,1;{\BR})\\
&\\
r_{2,1}^{X}\downarrow\quad&&r_{2,1}^{A}\downarrow\quad\\
&\\
H^{1,1}(X,{\BR}(1))&{\buildrel [w]^{\ast}\over\to}&
H^{1,1}(A,{\BR}(1))\\
\end{matrix}
\]
But $r_{2,1}^{X}$ and $[w]^{\ast}$ are surjective, whence
$r_{2,1}^{A}$ is surjective.
\end{proof}
 
\section{Proof of \thmref{THM101}}\label{SEC011}

The construction in \secref{SEC007}-\ref{SEC009}
dealt with $K_{1}$ classes on a general $K3$ surface $X$
that degenerate to $K_{0}$ classes on a special $K3$ surface
$X_{0}$ with maximum Picard number. Now suppose that we are
given $K3$ surfaces $X_{1},\ldots,X_{M}$, and $K_{1}$ classes $\xi_{j} =
\sum_{\alpha_{j}}\big(f_{\alpha_{j}},D_{\alpha_{j}}\big)$
with $\sum_{\alpha_{j}}\text{div}(f_{\alpha_{j}}) = 0$ on $X_{j}$, 
$j=1,\ldots,M$. Then we can form a $K_{1}$ cycle on the
product $X_{1}\times \cdots \times X_{M}$ by the prescription
$$
\xi = \sum_{\alpha_{1},\ldots,\alpha_{M}}\big(\text{Pr}_{1}^{\ast}
f_{\alpha_{1}}\cdots \text{Pr}_{M}^{\ast} f_{\alpha_{M}},
D_{\alpha_{1}}\times\cdots \times D_{\alpha_{M}}\big),
$$
where $\text{Pr}_{j}: D_{\alpha_{1}}\times\cdots \times D_{\alpha_{M}}
\to D_{\alpha_{j}}$ is the canonical projection. Note that
$\xi$ determines a class in $\text{CH}^{M+1}(X_{1}\times
\cdots \times X_M,1)$, and further note that as the $X_{j}$'s degenerate to
the special $X_{0}$, the class $\xi$ degenerates to an algebraic
cycle lying in $H^{1,1}(X_{0},{\BR}(1))^{\otimes M}$. Indeed
one can find $20^{M}$ such $\xi$'s, which degenerate to a basis
of $H^{1,1}(X_{0},{\BR}(1))^{\otimes M}$. Thus for general
$X_{1},\ldots,X_{M}$, $H^{1,1}(X_{1},{\BR}(1))\otimes\cdots
\otimes H^{1,1}(X_{M},{\BR}(1))$ lies in the image of
the real regulator $r_{M+1,1}: \text{CH}^{M+1}(X_{1}\times
\cdots\times X_{M},1;{\BR}) \to H_{\mathcal D}^{2M+1}(X_{1}\times
\cdots\times X_{M},{\BR}(M+1))$. Note that in particular,
in light of the previous section, this implies that for
a general product $E_{1}\times\cdots\times E_{2M}$ of
elliptic curves, $H^{1,1}(E_{1}\times E_{2},{\BR}(1))\otimes
\cdots \otimes H^{1,1}(E_{2M-1}\times E_{2M},{\BR}(1))$
lies in the image of the real regulator. Let
general mean as in \secref{SEC002},
with respect to the countable real
analytic Zariski topology. We are now in a position to prove
the following.

\bigskip

\begin{thm}\label{THM110}
Let $X = E_{1}\times \cdots \times E_{k}$ 
be a general product of $k$ elliptic curves, i.e. in the countable
real analytic Zariski topology, and
${\rm CH}_{\rm ind}^{k}(X,1;{\BQ})$ the
space of indecomposables. Then:

\medskip

{\rm (i)} {\rm Level}$\big({\rm CH}_{\rm ind}^{k}(X,1;{\BQ})\big) \geq k-2$.

\medskip

{\rm (ii)} In particular for $k\geq 3$, there are an 
uncountable number of indecomposables
in the kernel of the regulator 
$$
{\rm cl}_{k,1}: {\rm CH}_{\rm ind}^{k}(X,1;{\BQ})
\to H_{\mathcal D}^{2k-1}(X,{\BQ}(k)).
$$
\end{thm}

\begin{proof}
Using the notation just preceeding (11.1),
as well as the terminology in (5.5) and (5.7),
we put $k=M+1$, $X = E_{1}\times\cdots\times  E_{k}$, and let $S = 
E_{k+1}\times\cdots\times E_{2k-2}$. Note that
$\dim S = \ell -1$, where $\ell = k-1$.
Then the regulator map under consideration
is 
$\text{cl}_{k,1} : \text{CH}^{k}(X,1) \to H_{\mathcal D}^{2k-1}(X,{\BQ}(k))$.
Note that for a general product of elliptic curves, 
$$
N^{1}H^{k}(X,{\BQ}) = \sum_{j=1}^{k}\bigg[\frac{dx_{j}d\overline{x}_{j}
}{|y_{j}|^{2}}\biggr]\bigcup H^{k-2}(X,{\BQ}),
$$
where $E_{j}$ is defined in affine coordinates
by the equations $y_{j}^{2} = h_{j}(x_{j})$,
for general cubic polynomials $h_{j}(x_{j})$, and $j=1,\ldots,k$.
It is obvious that $H^{\{k,k-1,1\}}(X) \ne H_{N}^{1,k-1}(X)$.
Thus by \thmref{THM507}, 
\[
\text{Level}(\text{CH}_\text{ind}^{k}(X,1;{\BQ})) \geq \ell - 1 = k-2.
\]
Since the Hodge conjecture is known for products of elliptic
curves, it follows by 
\remref{REM508}
that $\ker \text{cl}_{k,1}$ contains an uncountable
number of indecomposable elements, provided that $k>2$.
\end{proof}

\appendix

\section{}

The analogue of a $K3$ surface in dimension $1$ is an elliptic
curve, and the group of interest is $K_{2}$. We prove
a version of the Hodge-${\mathcal D}$-conjecture for general
elliptic curves.
Let $X$ be a compact Riemann surface of genus $g$,
and let $f,g \in \BC(X)^\times$. For a real form
$\omega\in H^1(X,\BR)$, the integral
$$
\{f,g\} \mapsto \int_X\log|f|d\log|g|\wedge\omega
$$
induces a map on $K_2(X)$; more explicitly a map
\[
\begin{split}
r_{2,2} : {\rm CH}^2(X,2;\BR) := H_{\rm Zar}^0(X,{\mathcal K}_{2,X})
\otimes \BR
&\to H^1(X,\BR)^\vee\\
& \simeq H^1(X,\BR(1)) \simeq
H_{\mathcal D}^2(X,\BR(2)).
\end{split}
\]
Up to a multiplicative constant and real isomorphism,
$r_{2,2}$ is the real Beilinson regulator \cite{Lw1}. It
is well-known that $r_{2,2}$ is zero for general
curves of genus $g>1$ \cite{Co1}, is nontrivial for
the case of general elliptic curves ($g=1$) \cite{Blo2}, \cite{Co1},
and is trivially surjective for $g=0$. We sketch a proof of:

\bigskip

\renewcommand{\thethm}{}
\begin{thm}[Hodge-$\mathcal D$ for Elliptic Curves]
If $X$ is a general elliptic curve in the real
analytic Zariski topology, then $r_{2,2}$ is surjective.
\end{thm}

\bigskip

\begin{proof}
Let $X$ be an elliptic curve
given in affine coordinates by the equation
$y^2 = h(x)$, where $h(x)$ is a cubic polynomial with
distinct roots. A basis for $H^1(X,\BR)$ is given by
$$
\omega_1:= \frac{dx}{y} + \frac{d\overline{x}}{\overline{y}}\quad;
\quad \omega_2:= \sqrt{-1}\biggl(\frac{dx}{y} - 
\frac{d\overline{x}}{\overline{y}}\biggr).
$$
Next, we consider
$$
f_1 := y+x\sqrt{-1}\quad ;\quad f_2 = y+x\quad ;\quad
g_1 = g_2 = x.
$$
We claim that for general $X$,
\[
\det
\begin{bmatrix}
\int_X\log|f_1|d\log|g_1|\wedge\omega_1&
\int_X\log|f_1|d\log|g_1|\wedge\omega_2\\
&\\
\int_X\log|f_2|d\log|g_2|\wedge\omega_1&
\int_X\log|f_2|d\log|g_2|\wedge\omega_2
\end{bmatrix}
\ne 0.
\leqno{(1)}
\] 
Now let us first assume that $X$ is given for which (1) holds,
and note that the rational functions $f_1,f_2,g_1,g_2$
can each be expressed in the form $L_1/L_2$, where
$L_j$ are homogeneous linear polynomials in the homogeneous
coordinates of $\P^2$ (and where $X\subset \P^2$).
Since $X$ has a dense subset of torsion points $X_{\rm tor}$, 
and by Abel's theorem, one can find
$\tilde{L}_j$ ``close'' to $L_j$, $j=1,2$, such that
$\tilde{L}_j \cap X \subset X_{\rm tor}$. Thus
$\tilde{L}_1/\tilde{L}_2$ is ``close'' to $L_1/L_2$. Thus one
can find $\tilde{f}_1,\tilde{f}_2,\tilde{g}_1,\tilde{g}_2$
for which 
$$
\biggl\{\big|{\rm div}(\tilde{f}_1)\big|\bigcup \big|{\rm div}(\tilde{f}_2)\big|
\bigcup \big|{\rm div}(\tilde{g}_1)\big|\bigcup
\big|{\rm div}(\tilde{g}_2)\big|\biggr\}\subset X_{\rm tor},
$$
and that by continuity considerations
\[
\det
\begin{bmatrix}
\int_X\log|\tilde{f}_1|d\log|
\tilde{g}_1|\wedge\omega_1&
\int_X\log|\tilde{f}_1|d\log|\tilde{g}_1|\wedge\omega_2\\
&\\
\int_X\log|\tilde{f}_2|d\log|\tilde{g}_2|\wedge\omega_1&
\int_X\log|\tilde{f}_2|d\log|\tilde{g}_2|\wedge
\omega_2
\end{bmatrix}
\ne 0.\leqno{(2)}
\]
From the general mechanism in \cite{Blo2}, one can complete
$\{\tilde{f}_1,\tilde{g}_1\}$, $\{\tilde{f}_2,\tilde{g}_2\}$
to classes $\xi_1,\ \xi_2\in {\rm CH}^2(X,2)$, for which
\[
\det
\begin{bmatrix}
r_{2,2}(\xi_1)(\omega_1)&
r_{2,2}(\xi_1)(\omega_2)\\
&\\
r_{2,2}(\xi_2)(\omega_1)&
r_{2,2}(\xi_2)(\omega_2)\\
\end{bmatrix}
\ne 0.\leqno{(3)}
\]
Thus modulo the claim in (1), we are done.
\bigskip
We sketch a proof of the claim. 
With regard to a volume element $dV$:
$$
d\log|x| \wedge \omega_1 = \frac{1}{2}\biggl(\frac{1}{x\overline{y}} -
\frac{1}{\overline{x}y}\biggr)dx\wedge d\overline{x} =
\frac{{\rm Im}(\overline{x}y)}{|x|^2|y|^2}dV\leqno{(4)}
$$
$$
d\log|x| \wedge \omega_2 = -\frac{\sqrt{-1}}{ 2}\biggl(\frac{1}{x\overline{y}} 
+ \frac{1}{\overline{x}y}\biggr)dx\wedge d\overline{x} =
-\frac{{\rm Re}(\overline{x}y)}{|x|^2|y|^2}dV\leqno{(5)}
$$
Now let us degenerate $X$ to the rational elliptic curve
$X_0$ given by $y^2=x^3$. Note that $X_0$ is given parametrically
by $(x,y) = (z^2,z^3)$, $z\in \BC$. Thus
$\overline{x}y = |z|^4z$, and up to a real positive multiplicative
constant times the standard volume element on $\BC$,
which we will denote by $dV_0$, (4) and (5) become:
$$
d\log|x| \wedge \omega_1 = \frac{{\rm Im}(z)}{ |z|^4}dV_0\quad;
\quad d\log|x| \wedge \omega_2 = -\frac{{\rm Re}(z)}{ |z|^4}dV_0.\leqno{(6)}
$$
Let ${\bf H} = \{z\in \BC\ \big|\ {\rm Im}(z)\geq 0\}$ be the
upper half plane.
Now one has the following formal calculations after degenerating to $X_0$,
and using symmetry arguments:
\begin{align}
\tag{7}
&\quad \int_{X_0}\log|f_1|d\log|g_1|\wedge\omega_1\\
\notag
& = \int_\BC
\log|z^3+\sqrt{-1}z^2|\frac{{\rm Im}(z)}{|z|^4}dV_0
= \int_\BC
\log|z+\sqrt{-1}|\frac{{\rm Im}(z)}{|z|^4}dV_0
\\
\notag
&= \int_{\bf H}\log\biggl|\frac{z+\sqrt{-1}}{\overline{z}+\sqrt{-1}}
\biggr|\frac{{\rm Im}(z)}{|z|^4}dV_0 \mapsto +\infty,
\end{align}
using the fact
\[
\biggl|\frac{z+\sqrt{-1}}{
\overline{z}+\sqrt{-1}}\biggr| > 1 \Leftrightarrow\ {\rm Im}(z)>0.
\]

$$
\int_{X_0}\log|f_2|d\log|g_2|\wedge\omega_1 
= -\int_\BC\log|z+\sqrt{-1}|\frac{{\rm Re}(z)}{|z|^4} dV_0 = 0.
\leqno{(8)}
$$

$$
\int_{X_0}\log|f_2|d\log|g_2|\wedge\omega_1 
= \int_\BC
\log|z+1|\frac{{\rm Im}(z)}{|z|^4}dV_0 = 0.\leqno{(9)}
$$
For the final calculation, put $w = z\sqrt{-1}$,
and note that ${\rm Re}(z) = {\rm Im}(w)$, and
that $|z+1| = |w+\sqrt{-1}|$. Then
\begin{align}
\tag{10}
\int_{X_0}\log|f_2|d\log|g_2|\wedge\omega_2
&= -\int_\BC\log|z+1|\frac{{\rm Re}(z)}{|z|^4}dV_0\\
\notag
&= -\int_\BC\log|w+\sqrt{-1}|\frac{{\rm Im}(w)}{|w|^4}dV_0\\
\notag
&= -\int_{\bf H}\log\biggl|\frac{z+\sqrt{-1}}{\overline{z}+\sqrt{-1}}
\biggr|\frac{{\rm Im}(z)}{|z|^4}dV_0 \mapsto -\infty.
\end{align}
Notice that the singularities in the integrals in (7) and (10)
occur over the singular point $z=0$ of the singular curve
$X_0$, as expected. By using the Lebesgue theory of integration,
we can make the calculations in (7)--(10) above more precise.
First, by using
the projection $(x,y)\mapsto x$, we have a double covering
$X\to \P^1$. Thus for $f,g\in \BC(X)$, and $\omega = \omega_1$
or $\omega = \omega_2$, we can express $\int_X\log|f|d\log|g|\wedge
\omega$ as the integral of some Lebesgue integrable
function $H(x)$ over $\P^1$. Next, by converting to
polar coordinates, viz. $x = {\rm e}^{t\sqrt{-1}}$, we can Fubini
integrate in $t\in [0,2\pi]$ and $r \in [0,\infty]$. 
Let $h(r)$ be the result of integrating
$H(x)$ with respect to $t$ over $[0,2\pi]$. As $X$ degenerates
to $X_0$, we can construct a sequence $\{h_n(r)\}$ which
limits to $h_\infty(r)$ over $X_0$. In the cases of (7)--(10),
we have that $h_\infty(r)$ is either zero, nonnegative, or nonpositive. 
By using the the standard Lebesgue integral limit theorems,
we arrive at the
claim in (1), and hence the theorem above.
\end{proof}

\end{document}